\documentclass[12pt]{amsart}
\usepackage{latexsym,fancyhdr,amssymb,color,amsmath,amsthm,graphicx,listings,comment}
\usepackage[section]{placeins}
\let\paragraph\subsection
\definecolor{yellow1}{rgb}{1,1,0.8}
\def\satz#1{ \vspace{2mm} \begin{center} \fcolorbox{yellow1}{yellow1}{ \parbox{15.5cm}{{\bf Summary:} #1}} \vspace{2mm} \end{center} }
\title{Finite topologies for finite geometries}
\author{Oliver Knill}
\date{January 8, 2023}
\address{Department of Mathematics \\ Harvard University \\ Cambridge, MA, 02138 }


\begin{document}
\maketitle

\begin{abstract}
Without leaving finite mathematics and using finite topological spaces only, 
we give a definition of homeomorphisms of finite abstract simplicial complexes or finite graphs. 
Besides exploring the definition in various contexts, we add 
some remarks like that the general Lefschetz formula works for any continuous map
on any finite topological space. We also noted that any higher order Wu characteristic
as well as their cohomology are topological invariants which are not homotopy invariants. 
Energy theorems allow to express these topological invariants in terms of 
interaction energies of local open sets. 
\end{abstract}

\section{About}

\paragraph{}
When exploring finite geometries using finite topological spaces, one is challenged with the fact
that homeomorphic finite topological spaces have the same cardinality, which is too rigid.
A finite metric space produces the discrete topology which is too fine. 
For the geodesic metric on a graph for example, the topology is totally disconnected 
and so does not reflect at all the connectivity of the graph. As probably first realized
in 1937 by Alexandroff \cite{Alexandroff1937}, non-Hausdorff finite topologies still can capture essential 
parts of a topology that is usually only explored using geometric realizations. How can one avoid 
geometric realizations and still have a workable definition of homeomorphism?
We have explored a related notion in \cite{KnillTopology} using covers but use now
classical finite topologies, accepting the fact that all reasonable finite topological spaces
are Alexandroff and naturally non-Hausdorff if they capture connectivity properties of 
the space under consideration. 

\paragraph{}
This working document has grown a bit longer than anticipated but has been a seed for further results like
Green function formulas. As a remedy, we added summaries at the end of each section. 
The write-up is a contribution to a program of 
replacing continuum geometries in a finite set-up but with as little changes in notation as possible. 
We want to avoid geometric realizations because using the continuum is a rather
serious step. It is not just a philosophical obsession: the mathematics of topological
manifolds has told lessons like that there are finite geometries $G$ - and example is the simplicial 
complex obtained by taking the join of a discrete circle with a homology $3$-sphere - 
which produces a geometric realization which is classically homeomorphic to the standard 
$5$-sphere $H$ even so from any finite point of view they are not homeomorphic. They are
geometries which have homeomorphic geometric realizations but should not be considered
homeomorphics. In some sense the mathematics of the Hauptvermutung \cite{Hauptvermutung} has indicated that the
continuum can lead  to identifications which are not expeected. A finite topology more honestly preserves details which 
geometric realizations do not see any more. One could of course use piecewise linear geometry to capture
what finite topology does, but also from a computer science point of view, it is desirable to have
finite objects and finite data to deal with only. All geometric objects  are faithfully implemented using a finite
amount of information. The axiom of infinity is never used. 

\paragraph{}
Finite topologies by definition are always Alexandroff spaces \cite{Alexandroff1937}, 
meaning that every point $x$ has a smallest neighborhood $U(x)$. When working with a sheaf over
such a topology, we do not need to conceptualize direct limit constructions like
``stalks" or ``germs". In the case of simplicial complexes, the smallest atomic ``Planck units
of space" containing $x$ is known as the ``star" $U(x)$ of the simplex $x$. 
As connection calculus \cite{KnillEnergy2020} illustrates, the topology of intersections 
of stars $U(x) \cap U(y)$ can be complicated, even so each $U(x),U(y)$ is contractible. 
Indeed, the Euler characteristic of $U(x) \cap U(y)$ agrees up to a sign with the matrix 
entries of the inverse $g(x,y)$ of the connection matrix $L(x,y)$, which is $1$ if the 
simplices intersect and $0$ else. In analogy with the {\bf Green} functions in classical frame works, 
these numbers $g(x,y)$ must be thought of the potential energies between $x$ and $y$. They 
can be rather arbitrary for large dimensional spaces. One of the reasons why small dimensional 
topology is so much different from larger dimensional ones is that there can be surprises in
the topology of local ``atomic parts of space". While simple in small dimensions, the 
intersection $U(x) \cap U(y)$ can be entangled in a rather complicated way in higher dimensions. 

\begin{figure}[!htpb]
\scalebox{0.8}{\includegraphics{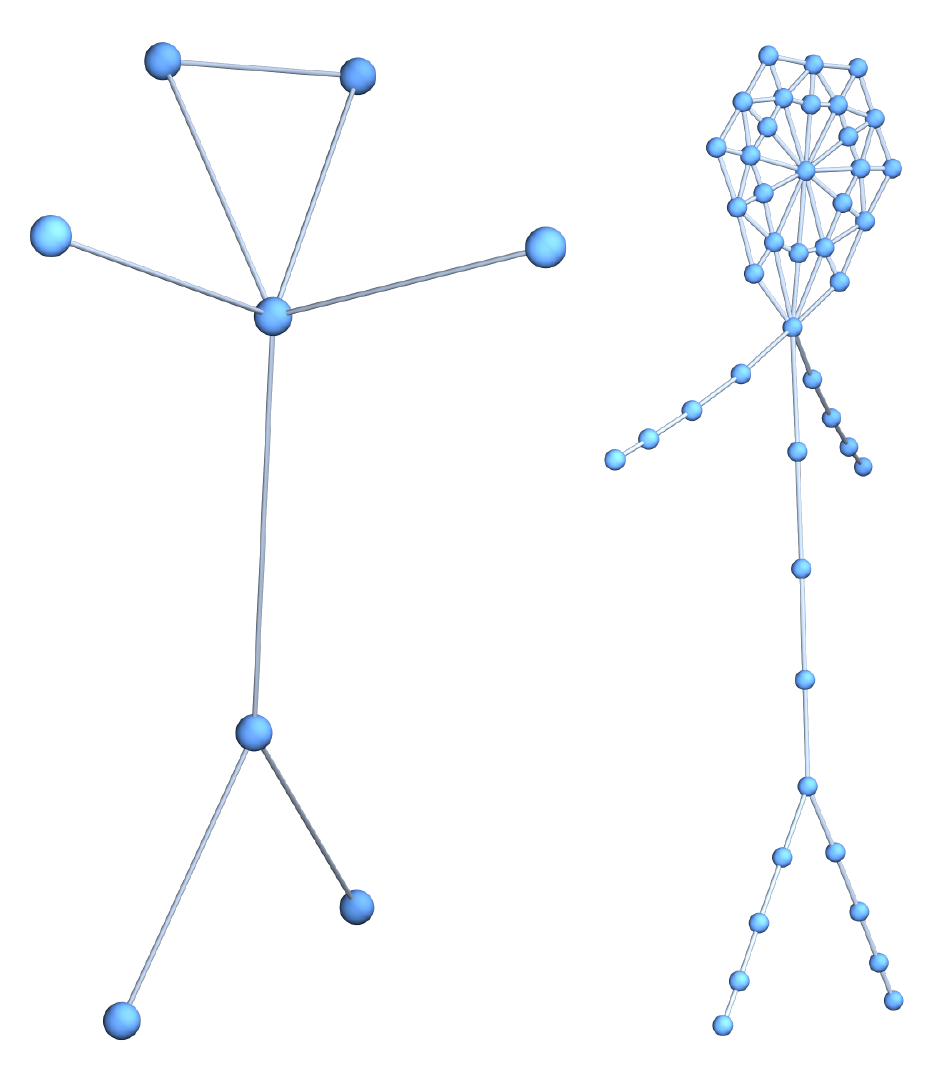}}
\caption{
The figure shows a finite simple graph $G$ and its second Barycentric refinement $G_2$. 
The Whitney simplicial complex $\mathcal{G}$ of $G$ has $17$ sets $x$, 
leading to $17$ basis elements $U(x)=U(x)$. These open sets in $\mathcal{G}$ are minimal 
and called the stars of $x$. The closure $\overline{U(x)}$ is the unit ball and its boundary 
$S(x) = B(x) \setminus U(x)$ is the unit sphere. 
The basis $\mathcal{B}$ generates a finite topology $\mathcal{O}$ with 3032 open sets 
and 3032 closed sets. Most of the $2^{17}=131072$ possible subsets of $\mathcal{G}$ 
are neither open nor closed. The topology is not Hausdorff: one can not separate 
points which intersect. As every finite topological space, it is
Alexandroff: every point $x$ has a smallest neighborhood $U(x)$. }
\end{figure}

\paragraph{}
Graphs and simplicial complexes and finite topological spaces
all provide model frame works in finite geometry. The categories are 
closely related: if one of them is given, one can construct
relatives in the other classes. One can get from a graph to a 
simplicial complex with the {\bf Whitney functor} by assigning to the graph 
the vertex sets of complete subgraphs. The \v{C}ech nerve construction
produces from a topological space a simplicial complex. 
\footnote{The Whitney complex is also known as the {\bf face complex}, clique complex, 
flag complex, or face poset.} 
From a simplicial complex, one can then construct a graph in which vertices are the sets of the 
complex and where two sets are connected if one is contained in the other. 
Switching forth and back between complexes and graphs produces a {\bf Barycentric refinement} 
of the topology. After identifying Barycentric refined complexes,
simplicial complexes or graphs can serve the same purpose.
We always construct the topology $\mathcal{O}$ on the simplicial complex $\mathcal{G}$ and denote
individual points in $\mathcal{G}$ with $x$. The sub-simplicial complexes of 
$\mathcal{G}$ are then the closed sets. In algebraic geometry, a similar 
constructions has led to the Zariski topology. 

\paragraph{}
While we deal here only with finite sets, most could be generalized to 
{\bf locally finite} complexes, meaning that there is an upper bound on the number
of elements in an atom $U(x)$ and again define the basis by the smallest neighborhoods $U(x)$.
This {\bf local finite assumption} corresponds in the continuum to the step to restrict to
{\bf paracompact} topological spaces, spaces where every open cover has a locally finite refinement. 
As for references in topology, see \cite{BourbakiTopologie,Hatcher,munkres,spanier} for topology, and 
especially \cite{alexandroff}, a text already using abstract simplicial complexes 
(introduced 1907 by Dehn and Heegaard \cite{BurdeZieschang}) and not the more commonly used definition 
using geometric realization. Veblen \cite{Veblen} in 1922 defined a {\bf neighborhood} of a k-simplex but
still used geometric realizations and cites Poincar\'e for introducing the notion of homeomorphism. 
Veblen also used already the terminology of ``stars". Every topology on a finite set is always an 
{\bf Alexandroff topology}. This is a topology, where points have smallest neighborhoods or 
alternatively where arbitrary intersections of sets are open too. The notion of ``Stars" was
established  in combinatorial topology like \cite{Alexandroff1937} but also entered some calculus 
textbooks like Whitney \cite{Whitney1957}. Still, as most texts, even Alexandroff
look at it primarily at geometric realizations of cell complexes. 
Alexandroff topology generalizes the {\bf co-finite topology} for $0$-dimensional complexes 
or the {\bf order topology} of a general partially ordered set with a basis 
$U(x)=\{ y, y \geq x \}$ \cite{Wachs2004}.

\paragraph{}
As for simplicial complexes, see \cite{FerrarioPiccinini} or \cite{MunkresAlgebraicTopology}.
Simplicial complexes appearing in graph theory are covered in \cite{JonssonSimplicial}. The most
important link between graphs and complexes is definitely the Whitney functor which assigns to 
a graph a complex which exactly has the topological properties which the graph suggests without
leading to ambiguities like what we consider to be a face.
In algebraic topology, the actual topology generated by the star basis has not obtained the 
attention it deserves but it appears, for example in \cite{spanier} on page 311.
The topology defined in \cite{AJK} is defined on the vertex set and of different nature, 
as even connectivity properties are different. See \cite{May2008} for a review on finite topological
spaces. Also \cite{Stallings1983} is completely unrelated because the 
star of a vertex is defined by Stallings as the set of edges attached to it and 
graphs considered one-dimensional objects, a common perspective in the 
20th century. There are quite a few other discrete frame works in finite geometry. We should mention
Ivashchenko \cite{I94,Evako1994} who translated Whitehead's homotopy notion into concrete procedures
in graph theory. It has been simplified in \cite{CYY} and crucial for defining what a 
``sphere" is combinatorially. The Morse approach is used in Forman's discrete 
Morse theory \cite{forman95,Forman1999}. In discrete combinatorics, one sometimes also looks a 
abstract simplicial complexes, in which the $\emptyset$ is included. 
We use the frame work, where $\emptyset$ is not considered to be a simplex but where it is 
considered to be a $(-1)$-dimensional sphere. The empty set itself is a simplicial complex,
as it fulfills the axiom, but it also does not contain the empty set. All these definitions 
are compatible with the continuum, where $d$-spheres also have Euler characteristic 
$1+(-1)^d$ even for $d=-1$ and where simplices all are contractible and have 
Euler characteristic $1$. 

\paragraph{}
We have explored the problem of defining homeomorphism within finite mathematics
for a few years already, the first time more seriously in 2014
\cite{KnillTopology}. Stars came up for us especially in the context of the 
{\bf Green star identity}, which explicitly gives the matrix entries
$g(x,y) = \omega(x) \omega(y) \chi(U(x) \cap U(y))$ of the inverse of the connection matrix 
$L(x,y) = \chi(\overline{ \{ x \cap y \} })$ attached to a simplicial complex $\mathcal{G}$ 
with simplices $x$ with $\omega(x) = (-1)^{{\rm dim}(x)}$ and Euler 
characteristic $\chi(A) = \sum_{x \in A} \omega(x)$ of a subset $A$ of $\mathcal{G}$. 
The {\bf connection matrix} $L$ involves the {\bf Euler characteristic} of closed sets
like $\overline{ \{x\} }$. Its inverse matrix $g$ involves the Euler characteristic 
of open sets $U(x)$. See \cite{AmazingWorld,GreenFunctionsEnergized,KnillEnergy2020}.
The sets $U(x)=W^+(x)$ are open sets and $\overline{\{x\}} = W^-(x)$ are closed sets. 
An important role plays the {\bf closure} $B(x) = \overline{U(x)}$ called the {\bf unit ball}.
Its {\bf boundary} $S(x) =\delta U(x) = \overline{U(x)} \setminus U(x)$ called
{\bf unit sphere}. They all are closed sets and so  simplicial complexes. 

\begin{figure}[!htpb]
\scalebox{0.7}{\includegraphics{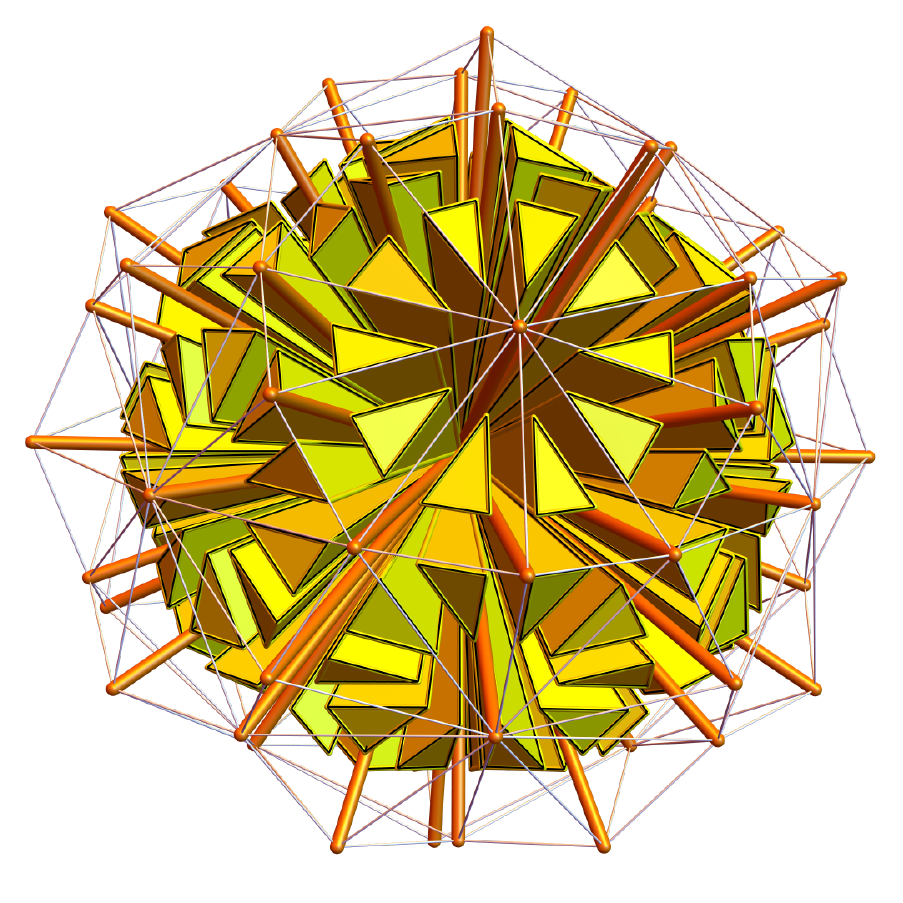}}
\scalebox{0.7}{\includegraphics{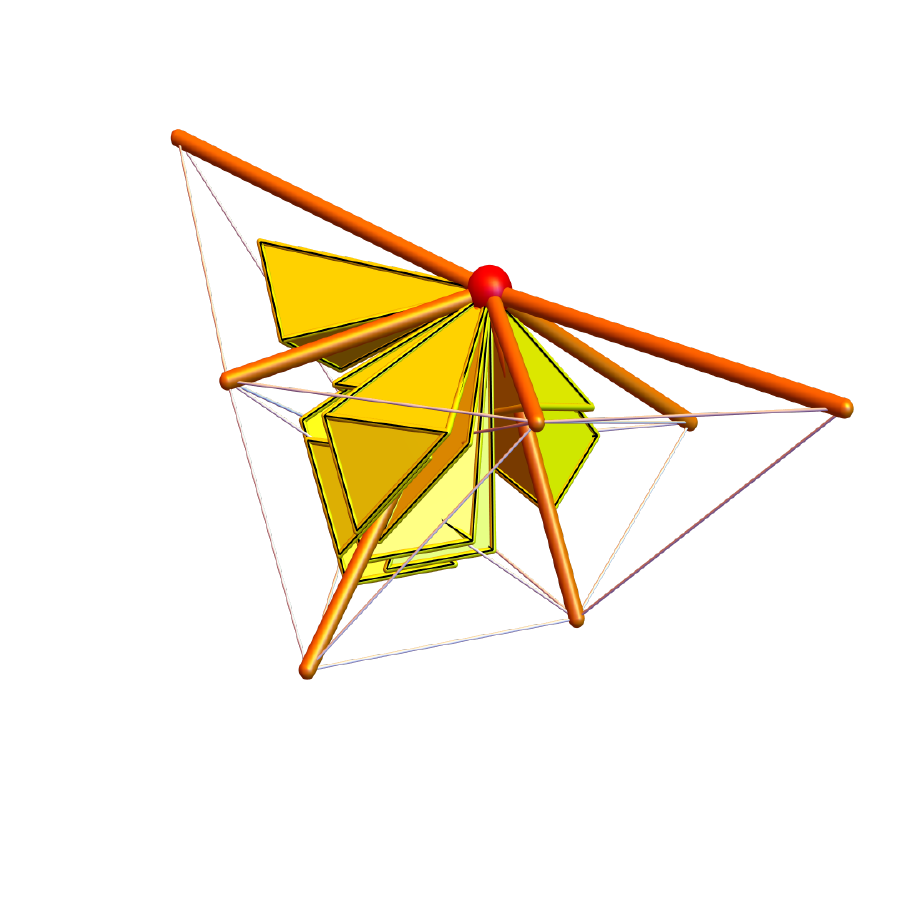}}
\caption{
A visualization of the star $U(x)$ in two examples. Unlike what the picture
suggests, we do not look at geometric realizations however. 
The star of a point $x \in \mathcal{G}$ consists of all the simplices
which contain $x$: that is $U(x) = \{y, x \subset y \}$. 
}
\end{figure}

\paragraph{}
In the case of a {\bf finite simple graph} $G$, we could also look at the 
topology $\mathcal{O}$ on its {\bf Whitney complex} $\mathcal{G}$, where the 
closed sets are the {\bf simplicial complexes coming from subgraphs} of $G$. 
This is a slightly rougher finite topology, because not all simplicial complexes 
are Whitney complexes of graphs. We are currently under the impression that the notion of 
homeomorphism discussed here is new. It has shifted quite a bit while we were writing this text. 
A useful definition needs to be simple and lead to all expected results. The goal
had been to get a definition of homeomorphism which generalizes the definition used 
for one-dimensional complexes in topological graph theory \cite{TuckerGross},
and which has the ability to identify different triangulations of obvious manifolds like
the icosahedron and octahedron. In one dimensions, the notion of homeomorphism is old and
enters for example the Kuratowski theorem:
two graphs $G,H$ are called {\bf graph homeomorphic} if there exists a graph isomorphism between 
some edge refined versions of $G$ and $H$. The story relating the finite and infinite 
is a bit tricky: more than a hundred years of work in the context of the {\bf Hauptvermutung} has lead 
to surprises in the relation between finite and infinite models. 
Finite topology can contain more information than 
the topology to geometric realizations. Using the topology from realizations 
in Euclidean spaces allows for some surprising homeomorphisms in higher dimensions. 

\paragraph{}
Having a topology on simplicial complexes is useful as it allows to reformulate classical 
results in a more familiar language but within a finite frame work.
For example: given any simplicial complex $\mathcal{G}$ and any 
continuous map $f: \mathcal{G} \to \mathcal{G}$, then the {\bf Lefschetz fixed point
formula} \cite{brouwergraph} $\sum_{x \in \mathcal{F}} i_f(x) = \chi_f(\mathcal{G})$ holds, where
$\mathcal{F}$ is the set of fixed points of $f$, where the index is defined as 
$i_f(x) = \omega(x) {\rm sign}(f|x)$ with ${\rm sign}(f|x)$ is the signature of the 
permutation which $f$ induces on the simplex $x$,
and where $\chi_f(\mathcal{G})$ is the super trace on cohomology, which is defined as
$\sum_{k \geq 0} (-1)^k {\rm tr}(U_f | {\rm ker}(L_k))$, where $L_k$ is the Hodge Laplacian $L=dd^*+d*d$ restricted to
the linear subspace of functions on $k$-dimensional simplices, where $df(x)= \sum_{y, |y|=|x|-1} {\rm sign}(y|x) f(y)$
is the {\bf exterior derivative} and $U_f g(x) = g(f(x))$ is the linear {\bf Koopman map} that $f$ 
induces on functions. 

\paragraph{}
Note that the Lefschetz fixed point theorem \cite{brouwergraph} holds for {\bf all simplicial complexes} 
and {\bf all continuous functions}. It generalized the theorem \cite{NowakowskiRival}, which is
a result in one dimensions. In the continuum, one needs assumptions, like that there are only 
finitely many fixed points. The discrete theorem had been formulated for graph endomorphisms 
\cite{brouwergraph} which produce continuous maps on the corresponding 
simplicial complex but the proof works also for continuous maps meaning for example that the map 
can contract. The Lefschetz fixed point theorem has two special cases: the first case is if $f$ is the 
identity, where it becomes the {\bf Euler-Poincar\'e formula}
$\sum_x \omega(x) = \sum_k (-1)^k b_k$, where $b_k$ are the Betti numbers, the dimensions of the kernels of $L_k$. 
An other special case is if the cohomology is trivial, meaning that only constant functions are
in the kernel of $L$. This applies for example if $\mathcal{G}$ is an arbitrary contractible complex and
among manifolds with boundaries if $\mathcal{G}$ is a $k$-ball. 
This is the {\bf discrete Brouwer fixed point theorem}: every continuous map 
on a finite abstract simplicial complex that is a $k$-ball has at least one fixed point. An other example
where we always have fixed points is if $\mathcal{G}$ is an even-dimensional sphere and if $f$ is continuous 
but preserves the orientation of the maximal simplices. These are results in finite mathematics. At no 
point, the concept of infinity is used. 

\paragraph{}
One of the main points was to have a clear definition what we mean with homeomorphism in finite 
topological spaces. Once one has such a notion, one can see what {\bf topological invariants are}. 
\cite{Bott52} defined {\bf combinatorial invariants} as properties invariant under Barycentric refinement.
We are especially interested in numerical quantities, that are not homotopy invariants. An example
in the continuum is the {\bf analytic torsion} \cite{KnillTorsion}, adapted from the continuum
\cite{RaySinger1971}. There are other properties we believe to be topological like being a Dehn-Sommerville
space \cite{DehnSommerville}. For all higher characteristics \cite{valuation}, 
which are only defined in the discrete so far, there are topological expressions which 
could be used to compute them in the continuum as we write down expressions which hold for 
the smallest open sets which exist in the topology. For the second characteristic, the {\bf  Wu characteristic} 
$\omega=\omega_2$, we know $omega(B)=(-1)^k$, where $k$ is the dimension of a ball $B$. One way to see this is that for
d-manifolds with boundary $\omega(B) = \chi(B) - \chi(\delta B)$ which in the case of a $d$-ball
is by the {\bf Euler gem formula} $1- (1+(-1)^{d-1}) =(-1)^d$. To compute the Wu characteristic of
a space, cover it by balls (which can have different dimensions but should match the dimension of 
the covered part. Within part of a ball not covered by different balls, the dimension of the underlying
space should be the same than the ball. What makes Wu characteristic compatible with topology is that
we have the {\bf valuation formula} $\omega(U \cup V) = \omega(U) + \omega(V) - \omega(U \cap V)$. 
This allows us to glue different parts together. Note that this valuation formula does not work for
closed sets which together with a formula $\sum_{x,y} \omega(x) \omega(y) \omega(U(x) \cap U(y))$ 
is the major reason why Wu characteristic is a topological invariant. 

\satz{A {\bf finite abstract simplicial complex} $\mathcal{G}$ is a finite set of
sets $x$ such that $\mathcal{G}$ is closed under the operation of taking 
finite subsets of $\mathcal{G}$. The cardinality $|x|$ of $x$
defines its {\bf dimension} ${\rm dim}(x)=|x|-1$. The {\bf star} $U(x)=\{ y, x \subset y\}$ of 
$x \in \mathcal{G}$ is the set of simplices containing $x$. It is declared to be {\bf open}.
The set $\mathcal{B}$ of stars $U(x)$ together with $\emptyset$ is a
{\bf topological basis} for a finite topology on $\mathcal{G}$. 
A complex $\mathcal{H}$ is called a {\bf d-ball}, if it is of the form 
$\mathcal{G} \setminus U(x)$, where $\mathcal{G}$ is a $d$-sphere. 
A complex $\mathcal{G}$ is {\bf contractible} if there exists $x \in \mathcal{G}$ 
such that both the {\bf unit sphere} $S(x)=\delta U(x)=\overline{U(x)} \setminus U(x)$ 
and $\mathcal{G} \setminus U(x)$ are contractible. An arbitrary set is 
called contractible if its closure is contractible. A complex is called a {\bf $d$-sphere} 
if it is a $d$-manifold and $\mathcal{G} \setminus U(x)$ is contractible for some $x$; a 
complex is a {\bf d-manifold} if every $S(x)$ is a $(d-1)$-sphere.
The complex $\mathcal{G}_1$ of all vertex sets of complete sub-graphs of the
graph $G_1=(V_1,E_1)$ with $V_1=\mathcal{G}$ and 
$E_1 = \{ (x,y), x \subset y, \; {\rm or} \;  y \subset x\}$ is called 
the {\bf Barycentric refinement} of $\mathcal{G}$. Define
$\mathcal{G}_n=(\mathcal{G}_{n-1})_1$. A simplex $x \in \mathcal{G}$ 
is {\bf locally maximal} if $x \subset y$ implies $y=x$. 
A complex $\mathcal{H}$ is declared to be a {\bf continuous image} of $\mathcal{G}$ 
if there exists a continuous surjective map $f: \mathcal{G}_n \to \mathcal{H}$ for some $n$
such that  (i) if $U(x) \subset \mathcal{H}$ for any locally maximal $k$-simplex 
$x \in \mathcal{H}$ has a pre-image whose closure is a $k$-ball and (ii) that every 
$f^{-1} S(x) \subset \mathcal{G}_n$ is homeomorphic to the unit sphere
$S(x) \subset \mathcal{H}$. Two complexes are {\bf homeomorphic} if each
is a continuous image of the other. 
These definitions are inductive either with respect to number of elements or dimension: 
the {\bf empty complex} = void  $0=\{ \}$ is the $(-1)$ sphere. The complex 
$1=\{ \{1 \} \}$ is contractible and the $0$-ball.}

\section{Topology}

\paragraph{}
An {\bf abstract finite simplicial complex} is a finite set $\mathcal{G}$ 
of non-empty sets that is closed under the operation of taking finite non-empty subsets.  
A {\bf finite simple graph} $G=(V,E)$ carries the {\bf Whitney complex}
$\mathcal{G}$ of $G$ which is the set of vertex sets of complete subgraphs of $G$. 
While not all simplicial complexes come from graphs in such a way,
\footnote{Examples are the $(k-1)$-dimensional {\bf skeleton complex} of the complete 
graph $K_k$ which is a $(k-1)$-sphere, also known as the boundary sphere of the simplex. } 
every simplicial complex $\mathcal{G}$ defines a {\bf finite simple graph} $G$ in which the sets 
$x$ of the complex are the vertices and where two sets are connected by an edge if one is 
contained in the other. If $\mathcal{G}$ came from a graph $G$ as a Whitney complex, 
the graph $G_1$ obtained from $\mathcal{G}$ is the {\bf Barycentric refinement} of $G$. 
Similarly, the Whitney complex $\mathcal{G}_1$ of $G_1$ is a simplicial complex, called the 
{\bf Barycentric refinement} of the simplicial complex $\mathcal{G}$. 
Since we can switch between graphs and complexes while doing Barycentric refinements,
the two concepts ``graphs" and ``complexes" can be interchanged. 
We like to keep both graphs and complexes and see them equipped with finite 
topological spaces. Simplicial complexes are attractive mathematical objects because they
have the {\bf simplest axiom system imaginable} in geometry: there is only one single axiom. 
Graphs on the other hand are unmatched in providing {\bf geometric intuition},
featuring {\bf accessibility}, and being supported by computer algebra systems,
much more than sets of sets.  
\footnote{Part of graph theory is accessible in secondary school education.
Simplicial complexes on the other hand tend to appear first in college topology 
or algebraic topology courses. The abstract version is more accessible because 
higher dimensional Euclidean spaces, usually only introduced in linear algebra courses,
are not invoked.}

\paragraph{}
As part of the definition, Barycentric refined objects are all homeomorphic so that  
from a topological point of view, they are identified. The goal is to use
standard notions of topology. We can not use the classical notion of homeomorphism for
finite topological spaces because this forces the finite topologies to be identical. 
We see all cyclic graphs to be homeomorphic for example or to see the icosahedron 
isomorphic to the octahedron as both are $2$-dimensional complexes which approximate 
under Barycentric refinements more and more spheres. 
We will call a complex $\mathcal{H}$ a {\bf continuous image} of $\mathcal{G}$ if there exists a 
continuous surjective map $f:\mathcal{G}_n \to \mathcal{H}$ such that for every $x$, 
the boundary $S(x)$ of $U(x)$ is homeomorphic to the boundary of $f^{-1}(U(x))$. This allows
to use induction with respect to dimension. We also require that for locally maximal simplices $x$
which have the property that the closure of the open set $f^{-1}(U(x))$ is a {\bf ball}, a 
simplicial complex which is obtained from a {\bf sphere} by removing an open set $U(z)$. 
If $\mathcal{G}$ is a continuous image of $\mathcal{H}$ and $\mathcal{G}$ 
is a continuous image of $\mathcal{H}$, the two spaces are considered homeomorphic. 

\paragraph{}
A finite abstract simplicial complex $\mathcal{G}$ is so always equipped with a 
finite topology $\mathcal{O}$ on $\mathcal{G}$. This is understood in the classical sense: 
a topology contains the empty set and $\mathcal{G}$, it is closed under finite intersections and
closed under arbitrary unions. In the finite case, we of courses can avoid the ``finite" word
but all we do here can be generalized to infinite but locally finite simplicial complexes. 
We stick with the finite, because the text should be seen as part of a larger and 
more ambitious project investigating the question: {\it which parts of geometry can be 
replaced with finite combinatorial notions?} We hope to be able to define
within finite sets of sets whether two simplicial complexes are homeomorphic or not and
point out that this is sharp than the softer equivalence relation given by 
homeomorphic geometric realizations. The classical ``homeomorphic notion" is 
too rigid for finite topological spaces as it forces a bijection between the atoms $U(x)$ making
up the basis. Applied to graphs it would require the graphs to be isomorphic.
We want cyclic graphs $C_n$ with $n \geq 4$ to be all 
homeomorphic for example. We want an edge refinement of a graph to be 
homeomorphic deformations and capture the notion of homeomorphism which is used in graph 
theory when graphs are considered one-dimensional simplicial complexes. 

\paragraph{}
If $\mathcal{G}$ is a finite abstract simplicial complex and $x \in \mathcal{G}$ is given, 
it defines the {\bf star} $U(x) = \{ w \in \mathcal{G}, v \subset w \}$.  
The collection $\mathcal{B}$ of all these stars together with the empty set
$\emptyset$ is a set of sets $\mathcal{B}$ that covers $\mathcal{G}$.
The collection $\mathcal{B}$ is also closed under intersections
because $U(x) \cap U(y) = U(x \cap y)$ if $x \cap y$ is not empty and $U(x) \cap U(y) = \emptyset$
else. It therefore defines a {\bf base} for a {\bf topology} $\mathcal{O}$ on $\mathcal{G}$.
By definition in point set topology, a {\bf topology} $\mathcal{O}$ is a set of
subsets of $\mathcal{G}$ which (i) contains $\emptyset$, (ii) contains $\mathcal{G}$ and which is
(iii) closed under finite intersections and (iv) closed under arbitrary unions.
A {\bf sub-base} of the topology is the set of sets $U(\{v\})$, with $v \in V=\bigcup_x x$. 
It is a sub-base because every base element $U(x)$ is an intersection of such sets
$U(x) = \bigcap_{v \in x} U(\{v\})$ and so generates the base from intersections.

\paragraph{}
A {\bf suspension} of graph $G$ is the {\bf Zykov join} of $G$ with the zero sphere $S_0$ (the graph with two
vertices and no edges). Doing this twice is a {\bf double suspension}. 
A {\bf rational homology 3-sphere} is a 3-manifold that has the same 
cohomology than a 3-sphere.  
\footnote{There are implementations of the homology 3-sphere with $16$ maximal simplices, 
leading to a complex with 392 simplices. The corresponding graph has 2552 edges. 
The double suspension has 394 vertices.}
The double suspension $G$ of a rational homology $3$-sphere is a concrete example of a
simplicial complex that is not a discrete $5$-sphere because not all unit spheres are 
spheres. But $G$ has a geometric realization that 
is a $5$-sphere by the {\bf double suspension theorem} of Edwards and Cannon. 
The finite topology can distinguish complexes in the discrete, 
which are indistinguishable using the tool of topological realizations. 
\footnote{PL geometry in the continuum would capture the finite topology too but it 
also would use the {\bf continuum} and as Euclidean spaces are used the concept of {\bf infinity}.} 
Edward \cite{Edwards1970} works with the {\bf Mazur homology 3-sphere} and shows that the 
double suspension is $S^5$. A consequence of the double suspension theorem is the 
existence of {\bf ``exotic triangulations"}: there are topological manifolds 
which are not equivalent to a piecewise linearly homogeneous polyhedron $P$
meaning that for any $x,y \in P$, there exists a piecewise linear 
homeomorphisms $h: P \to P$ such that $h(x)=y$. 
A triangulated topological manifold $M$ on the other hand only requires 
$h(x)=y$ which is a local homeomorphisms. Every smooth manifold 
has a PL structure. But a general topological manifold does not
need to be homeomorphic to polyhedra: Casson gave 4-manifold counterexamples.
Examples in higher dimensions have appeared more recently \cite{Manolescu}.

\paragraph{}
We start by defining some subsets in a finite abstract simplicial complex $\mathcal{G}$. 
Every $x \in \mathcal{G}$ defines the {\bf star} $U(x):=\{ y \in \mathcal{G}, x \subset y \}$
which is an open set and the {\bf core} $K(x)=W^-(x):=\{y \in \mathcal{G}, y \subset x\}$ which, 
unlike $U(x)$ in general, is always a sub-simplicial complex of $\mathcal{G}$ and so a {\bf closed set}.
The closure $B(x)$ of $\overline{U(x)}$ of $U(x)$ contains $\overline{\{x\}}=W^-(x)$ and is 
called the {\bf unit ball} of a point $x \in \mathcal{G}$. 
Its boundary $S(x) =\delta B(x) = \overline{B(x)} \setminus U(x)$ 
is a closed set called the {\bf unit sphere} of $x$. 
By definition, this is a closed set. In general, the boundary 
of any open set $\delta{U} = \overline{U} \setminus  U = \overline{U} \cap U^c$ is 
closed because it is the intersection of two closed sets. 
\footnote{The graph that can be constructed from the complex $S(x)$ 
is isomorphic to the subgraph of all points in distance $1$ to the vertex $x$ 
in the Barycentric graph $G_1$, the graph which is constructed from the complex 
$\mathcal{G}$. }

\paragraph{}
The unit sphere $S(x)=\delta U(x)$ is in the language of simplicial complexes, 
also known as the {\bf link} of $x$, (but it is usually only defined 
for $0$-dimensional $x$ so that we avoid the term). The open set 
$U(x)=W^+(x) = \{ y, x \subset y\}$ and closed set $K(x)=W^-(x) =\overline{ \{x\}} = \{ y, y \subset x\}$ 
are somehow {\bf dual} to each other. The closed set $W^-(x)$ is contained within $S(x)$. 
In the case when $x$ is {\bf locally maximal} meaning that it is not contained in a
strictly larger simplex, then $S(x)$ is the boundary complex of the simplex $x$ 
and so a sphere. A general complex $\mathcal{G}$ is a sphere:
there exists $y \in \mathcal{G}$ such that $\mathcal{G} \setminus U(y)$ is contractible
and also for every $y \in S(x)$, the unit sphere $S(y)$ within $S(x)$ is a 
co-dimension-one sphere again.

\paragraph{}
The terminology for graphs is similar.
In a finite simple graph $G$ and a vertex $v$, we call the graph generated by the 
vertices $w$ adjacent to $v$ the {\bf unit sphere} of $v$.
If $G$ comes from a simplicial complex $\mathcal{G}$, then each vertex has a
dimension and adjacent vertices are ordered. The sphere $S(x)$ 
is now the join of $S^+(v)$ generated by the vertices $w$ for which $w$ contains $v$ in
$\mathcal{G}$ and $S^-(v)$ is generated by all $w$ which are subsets of $v$. 
The sphere $S(v)$ of the graph is then the {\bf Zykov join} \cite{Zykov} of $S^-(x)$ and $S^+(x)$
because the vertex sets are the disjoint union and every element in $S^-(x)$ 
is connected to every element in $S^+(x)$. We like to think also of $S^{\pm}(x)$ 
as the {\bf stable and unstable manifolds} in the ``hyperbolic structure" defined by the 
{\bf Morse function} $f(x) = {\rm dim}(x)$ of {\bf Morse index} ${\rm dim}(S^-(x))+1$. 
Let us define this more generally. The join can also be defined directly for simplicial 
complexes by $\mathcal{G}+\mathcal{H} = \mathcal{G} \cup \mathcal{H} \cup 
\{ x+y, x \in \mathcal{G}, y \in \mathcal{H}\}$.

\paragraph{}
Let $R$ be an ordered ring like $\mathbb{Z},\mathcal{Q},\mathcal{Z}$.
A function $f: \mathcal{G} \to R$ is called {\bf locally injective} if $f(x) \neq f(y)$ 
for every $y \in S(x)$. A {\bf Morse function} on a complex $\mathcal{G}$ is defined
as a locally injective function function $\mathcal{G}$ which 
has the property that $S_f^-(x) = \{ y \in S(x), f(y)<f(x) \}$ is a $(k-1)$-sphere. 
Its {\bf Morse index} is $k$. If we do not want to refer to the graph and so to 
the Barycentric refined topology, we would require that $S_f^-(x)$ is
a simplicial complex which is a sphere.  
The function $f(x) = {\rm dim}(x)$ is a special Morse function. 
If $\mathcal{G}$ is a $d$-manifold, then the graph 
$S(x)$ in the graph $G_1$ is a $(d-1)$-sphere and agrees with the 
topological join of the two spheres $S^{\pm}(x)$. 
The geometric realization of a Zykov join of two graphs agrees with the 
{\bf topological join} of the geometric realizations

\satz{A finite abstract simplicial complex $\mathcal{G}$ carries a finite topology $\mathcal{O}$
in which the stars $\mathcal{B}$ form a basis. 
A finite simple graph carries so a natural topology on its Whitney complex. 
The topology is finite and so Alexandroff:
every point $x$ has a smallest neighborhood $U(x)$, the star of the simplex $x$. We like to 
think of them as atoms of space. The closure $B(x)$ of a star $U(x)$ is is called the 
unit ball. Its boundary $S(x)$ is called unit sphere of $x$.}

\section{Continuity}

\paragraph{}
The classical definition of continuity
can be applied immediately to functions $f$ between simplicial complexes 
$f:\mathcal{G} \to \mathcal{H}$ if we just silently assume the topology generated by stars
as the natural topology on the complex.  If we talk about continuous maps $f:G \to H$ of
graphs $G,H$, then rather looking at maps on the vertex sets $V(G),V(H)$, we look at maps
from its simplicial Whitney complex $\mathcal{G}$ of $G$ to the simplicial Whitney complex 
$\mathcal{H}$ of $H$. As usual in point-set topology, a map $f$ is called {\bf continuous} 
if the inverse $f^{-1}(A)$ of an open set $A$ in $\mathcal{H}$ is an open set in $\mathcal{G}$.

\paragraph{}
A map is continuous if and only if the inverse image of closed sets is closed. 
The standard definition of continuity works well, but it was necessary to modify the notion of 
``homeomorphismI" in the finite as classically, homeomorphic finite topological 
spaces are identical. It would be unacceptable for example to consider an icosahedron and octahedron 
as being topologically different, or to consider a cyclic graph with $6$ elements to be topologically 
different than a cyclic graph with $5$ elements. Even Barycentric refined complexes would not 
be homeomorphic with the narrow definition from point set topology, requiring the two continuous
maps which are inverses of each other. 

\paragraph{}
In the context of simplicial complexes, a continuous map is a bit more general than 
a {\bf simplicial map}. The later is a map from $\mathcal{G}$ to $\mathcal{H}$ that
preserves order and does not increase dimension if $x \subset y$ then $f(x) \subset f(y)$
and therefore satisfies ${\rm dim}(f(x)) \leq {\rm dim}(x)$. 
A simplicial map must map $0$-dimensional simplices to $0$-dimensional simplices. A 
continuous map $f: \mathcal{G} \to \mathcal{H}$ does not need to do that. A constant
map which has as an image a positive dimensional simplex is continuous but not a simplicial map. 
It does not necessarily map simplicial complexes into simplicial complexes. 
Most permutations $f: \mathcal{G} \mathcal{G}$ of a simplicial complex are not continuous. 
They scramble around the simplices without preserving the order. 
But simplicial maps are always continuous: 
if $f$ is a simplicial map, then the inverse image of any simplicial complex 
$\overline{\{y\})}$ consists of unions of simplicial complexes $\overline{\{x_k\}}$, 
which is closed. Because the inverse image of any closed set is a closed set, the map is continuous. 

\paragraph{}
Any map from a $0$-dimensional complex $\mathcal{G}$ to a complex $\mathcal{H}$ is always continuous 
because every set in $\mathcal{G}$ is both open and closed. 
Such a map neither does have to be injective, nor does it have to be surjective. 
An other extreme case is a constant map $f(x)=c$ from a complex $\mathcal{G}$ to a complex $\mathcal{H}$.
It is always continuous and a simplicial map if $c$ is zero dimensional: 
the set $f^{-1}(A)$ is either empty (if $c \notin A$) or then the entire space $\mathcal{G}$. 
The image $\{ c \}$ is however not open. So, even for finite topologies, a continuous map does not need to 
be an {\bf open map}, a map that transports open sets into open sets.

\paragraph{}
Sometimes it is good to look at maps defined within graphs $G=(V,E),H=(W,F)$ alone and not directly look 
at the simplicial complex. A map $f: V(G) \to V(H)$ is {\bf continuous graph map} $f:G \to H$ if 
$e=(a,b)$ in $G$ then either $(f(a),f(b))$ in $E(H)$ or then that $f$ collapses $e$ collapses to 
vertex $a=(a,a)$ in $V(H)$. Such a map $f$ lifts to a continuous map on the corresponding 
Whitney simplicial complexes $\mathcal{G} \to \mathcal{H}$. 
It actually even lifts to a simplicial map because zero dimensional parts get mapped into zero dimensional parts. 
Every continuous map between graphs as just defined leads so to a continuous map on the Whitney complex: 
given $x \in \mathcal{G}(G)$, define $f(x) = \bigcup_{v \in x} f(v)$. This is an element in the complex
$\mathcal{H}$ of $H$. If it was not, then there would exist $a,b \in x$ such that 
$(f(a),f(b)) \notin E$ nor $f(a)=f(b)$, contradicting the assumption that $f$ 
is a continuous graph map. Again, we should tell that there are continuous maps between simplicial 
complexes $\mathcal{G}(G), \mathcal{H}(H)$ of graphs $G,H$ which do not come from continuous 
graphs maps, the constant map to a positive
dimensional simplex is an example of a continuous map that does not come from a graph map
because a continuous graph map necessarily maps zero dimensional parts to zero dimensional parts. 

\paragraph{}
For a continuous map between graphs $f:G \to H$, the maximal dimension of the image 
graph $f(G) \subset H$ is always smaller or equal than the maximal dimension of $G$.
More generally, any simplicial map from a simplicial complex to an other simplicial complex
does not increase dimension on the image: 
if $f: \mathcal{G} \to \mathcal{H}$ is a continuous map between simplicial complexes then 
${\rm dim}(f(x)) \leq {\rm dim}(x)$. 
The dimension of the image of a continuous map can be strictly smaller of course: the constant map 
mapping every simplex in $\mathcal{G}$ to a single fixed vertex $v$ (zero-dimensional simplex)
in $\mathcal{H}$ is continuous and a simplicial map. The constant map to a positive dimensional 
simplex however is continuous but not a simplicial map. 

\paragraph{}
A {\bf graph homomorphism} \footnote{To have ``homomorphism" and ``homeomorphism" so 
close in the landscape of words is unfortunate but it is very much entrenched. 
The two terms rarely appear in the same context, but here they do.}
is a map between graphs $G=(V,E)$ and $H=(W,F)$ such that it maps $V$ into $W$ and $E$ into $F$. 
A graph homomorphism therefore maps simplices into simplices on the simplicial complex level 
and so defines a simplicial map and therefore a continuous map on its Whitney simplicial complexes:
the inverse image of an open set $U(x)$ in $\mathcal{G}(H)$ is open in $\mathcal{G}(G)$ 
because it is a disjoint union of sets $U(y_j)$
where $y_j$ are the set of simplices which are mapped into $x$. 
On the other hand, not every continuous map is a graph homomorphism because a graph 
homomorphisms by definition is not allowed to collapse an edge to a vertex. 

\satz{Continuity of maps between simplicial complexes 
is defined as usual in topology: the inverse image of an open 
set is open. For simplicial complexes, simplicial maps are continuous but the converse
is not necessarily true. There are continuous maps between simplicial complexes
which are not simplicial maps.
Graph homomorphisms define simplicial maps on their complexes and so are continuous too, 
but also here, the converse is not always true. 
In general, simplicial maps  or graph maps only can lower the maximal dimension: 
the dimension of $\mathcal{G}$ is larger or equal than the dimension of 
$f(\mathcal{G}) \subset \mathcal{H}$. But continuous maps between simplicial complexes do not
need to lower the dimension as the constant map $f(x)=c$ from a zero dimensional complex 
$\mathcal{G}$ to a positive dimensional complex $\mathcal{H}$ with a fixed positive 
dimensional $c \in \mathcal{H}$ shows. }

\section{Homeomorphism}

\begin{figure}[!htpb]
\scalebox{0.8}{\includegraphics{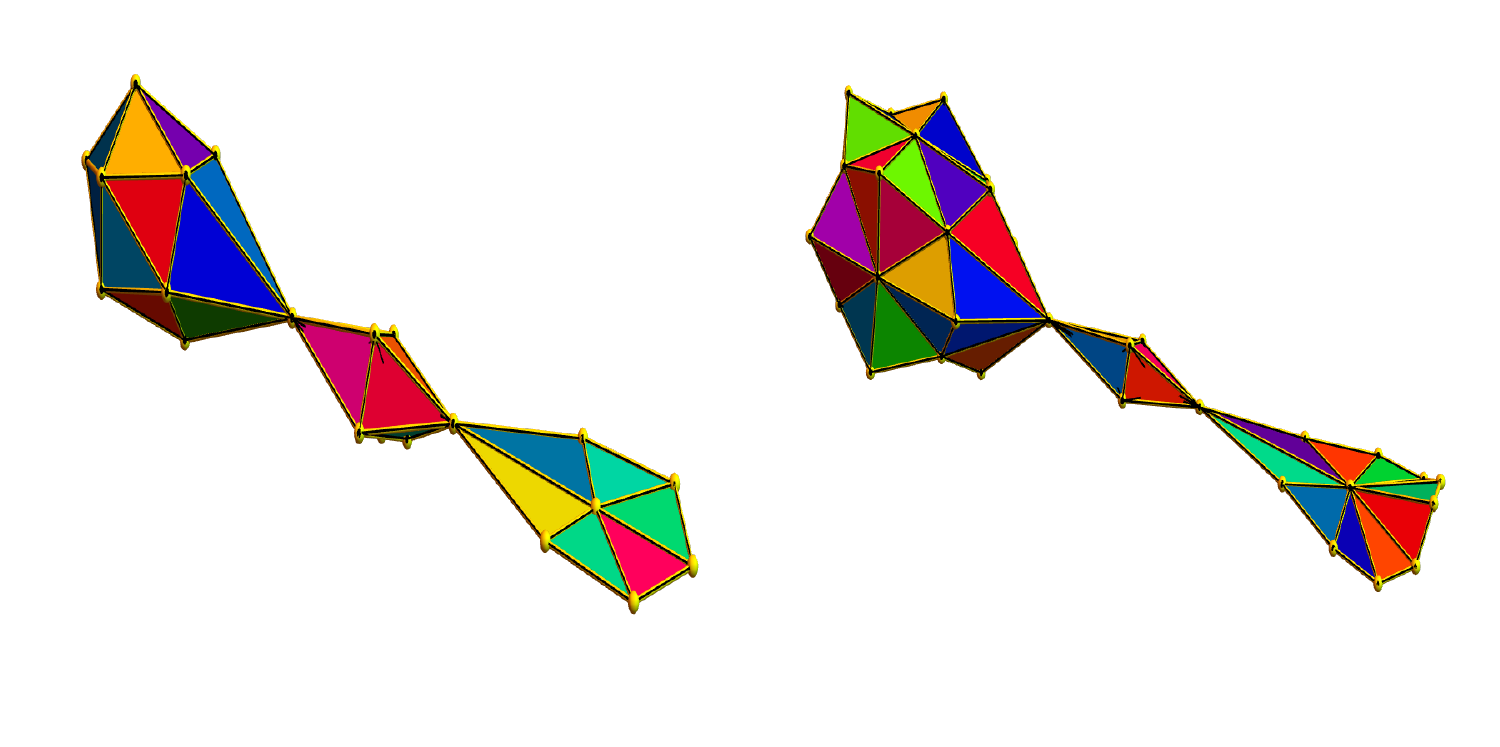}}
\caption{
Two homeomorphic finite abstract simplicial complexes $\mathcal{G},\mathcal{H}$ 
are displayed, with vertices, edges and triangles filled out. 
Both are a wedge sum of a 2-sphere, a 3-sphere and a 2-ball. In order to 
show that the two spaces are homeomorphic, one can first show that all
d-spheres are homeomorphic and all $d$-balls are homeomorphic and that if
two pointed spaces are homeomorphic then their wedge sums are homeomorphic. 
The wedge sum of two path graphs can be a star graph or a path graph and they are
not homeomorphic. One-dimensional complexes are homeomorphic if and only if 
they are classically homeomorphic, that is if one can get from one to the other
by a sequence of edge refinements or edge collapses that come from edge refinements. }
\end{figure}

\paragraph{}
As pointed out in the introduction, the classical notion of homeomorphism is too
rigid for finite topologies. Already pioneers like
Poincar\'e used in a combinatorial set-up at equivalence classes of 
finite geometries and considered Barycentric refinement complexes equivalent. 
\footnote{As pointed out earlier, PL-geometry \cite{RourkeSanderson} would do the job 
in the continuum. But we do not want to use Euclidean spaces, nor use of infinity. }
So, to start with, we assume that two geometries which are Barycentric refinement
of the other are homeomorphic. 
The {\bf Barycentric refinement} $G_1=(V_1,E_1)$ of a graph $G=G_0$ is a new
finite simple graph defined as follows: the vertex set of $G_1$ is $V_1=\mathcal{G}$. 
The Barycentric refinement of a simplicial complex $\mathcal{G}$ 
is the Whitney complex of the graph defined by $\mathcal{G}$. 
The edge set is the set of pairs $(x,y)$ for which either $x \subset y$ or $y \subset x$. 
We can iterate the Barycentric construction and look at the Barycentric refinements
$G_n = (G_{n-1})_1$ for every $n \geq 0$. We have just seen that a continuous map 
$G \to H$ can be lifted to a continuous map $G_1 \to H_1$. This means now that every
continuous map $G \to H$ between two graphs be lifted to a continuous map on the
Barycentric refinements $G_n \to H_n$. 

\paragraph{}
We repeat the new definition: if there is a continuous map $f: \mathcal{G}_n \to \mathcal{H}$ 
such that for every $y \in \mathcal{H}$, the boundary of $f^{-1}(U(y))$ is homeomorphic to 
the boundary $S(y)$ of $U(y)$ and for all locally maximal simplices $y$, the ball 
$B(y) = \overline{U(y)}$ in $\mathcal{H}$ has a pre-image $f^{-1}(B(y))$ which is a ball,
we say $\mathcal{H}$ is a {\bf continuous image} of $\mathcal{G}$. 
We say $\mathcal{G}$ and $\mathcal{H}$ are {\bf  homeomorphic} if 
$\mathcal{H}$ is a continuous image of $\mathcal{G}$ and $\mathcal{G}$ is a 
continuous image of $\mathcal{H}$.  
A {\bf d-ball} is a complex that is of the form $\mathcal{S}-U(z)$, 
where $\mathcal{S}$ is a d-sphere and $z \in \mathcal{S}$. 
Remember that a $d$-sphere $\mathcal{S}$ is a $d$-manifold which when 
punctured ($\mathcal{S} \setminus U(z)$) becomes contractible and that 
a $d$-manifold is a complex for which all unit-sphere is a $(d-1)$-sphere.
The empty complex is the $(-1)$-sphere and the 1 point complex contractible. 
\footnote{This is a {\bf recursive definition} as it refer to homeomorphism smaller dimensions.
In zero dimensions, homeomorphic means equal cardinality.}

\paragraph{}
An immediate consequence of the definition is that $0$-dimensional complexes are 
homeomorphic if and only if they have the same number of elements. The reason is that we require the
inverse image of every $B(x) = \{x\}$ to be a $0$-ball which is a one-point complex $K_1$. This
implies that $f$ must be injective. The map $f$ therefore has to be a bijection. 
This defines now both an {\bf equivalence relation} between simplicial complexes as well as 
for graphs: it is reflexive and symmetric. To see transitivity of the relation, 
note that we have now a chain of maps
$\mathcal{G}_n \to \mathcal{H}_m \to \mathcal{K}_l$ showing that $\mathcal{K}$ is a
continuous image of $\mathcal{G}$. Since the reverse holds also, the complexes $\mathcal{K}$ 
and $\mathcal{G}$ are homeomorphic, if $\mathcal{G}$ and $\mathcal{H}$ are homeomorphic
and $\mathcal{H}$ and $\mathcal{K}$ are homeomorphic. 

\paragraph{}
The notion of homeomorphism goes over to graphs $G$ if we look at the Whitney simplicial 
complex $\mathcal{G}$ attached to the graph $G$. The topology of $\mathcal{G}$ is then taken. 
By definition then, the Barycentric refinements $G_n$ of a graph $G$ 
are all homeomorphic to each other. We could also restate this by noting 
that $G$ and $H$ are homeomorphic if and only if $G_n$ and $H_n$ are homeomorphic 
for some $n$. As an example: two cyclic graphs $G=C_n$ and $H=C_m$ with $n,m \geq 4$ 
are homeomorphic. If $G$ is a triangulation of a compact manifold $M$ and $H$ is a triangulation
of a compact manifold $N$ such that $G,H$ have a common Barycentric refinement or common edge refinement,
then $G,H$ are homeomorphic and $M,N$ are combinatorially equivalent.
Any {\bf combinatorial invariants}, a term coined in \cite{Bott52} (meaning a property that
does not change under Barycentric refinements) must also be a {\bf topological invariant} 
(meaning a property that is the same for homeomorphic objects). 

\paragraph{}
It is useful to reformulate the notion of homeomorphism as the property that 
there exists maps $f: G_n \to H_m$ and $g:  H_m \to G$ which are both continuous
and such that the homeomorphism works also locally in that the smallest spheres have
pre-images which are homeomorphic and that the inverse of locally maximal unit balls
are actual balls. We still wonder whether the assumption on the locally maximal
balls can be avoided. The definition includes it so that we can prove things, 
like that if two complexes are homeomorphic and one of them is a manifold, 
then the other must  also be a manifold. 

\paragraph{}
Let us dwell on this a bit more:
if a map $f: \mathcal{G} \to \mathcal{H}$ is continuous, we could try to ask that
all unit spheres $S(x) = \delta U(x)$ are homeomorphic to $\delta f^{-1} U(x)$ and
use this alone as a recursive definition for homeomorphism. 
We would then start with the assumption that zero-dimensional complexes are 
homeomorphic if they have the same cardinality,
the property $\delta f^{-1} U(x)$ homeomorphic $\delta U(x)$ would then 
recursively lift the notion of homeomorphisms dimension by dimension. It is still
not clear whether this alternative definition alone would work.

\paragraph{}
While the definition of homeomorphism is by design symmetric, we could also explore and
try assuming only direction: checking this for $f: \mathcal{G}_n \to \mathcal{H}$ in one direction 
only might already determine that $\mathcal{G}$ and $\mathcal{H}$ 
are homeomorphic. This works if $\mathcal{G}$ and $\mathcal{H}$ are 
one-dimensional. In that case, we need a map from $\mathcal{G}_n$ to $\mathcal{H}$ 
such that the vertex degrees of $S(x) = \delta U(x)$ is the same than the vertex 
degrees of $\delta f^{-1} U(x)$ for every $x$. We come back to this question again at the end. 

\satz{The complex $\mathcal{H}$ is declared to be a continuous image of $\mathcal{G}$, if the natural
surjective map from some Barycentric refinement $\mathcal{G}_n \to \mathcal{G}$ factors
as $\mathcal{G}_n \to \mathcal{H}_m \to \mathcal{G}$ and the map induces
homeomorphisms on unit spheres and maximal unit balls have pre-images that are 
actual balls, punctured spheres. If two complexes are continuous images of each other, we call them 
homeomorphic. This notion defines an equivalence relation on complexes 
as well as an equivalence relation on graphs. }

\section{Closed}

\paragraph{}
Closed sets of a graph $G$ are subsets of $\mathcal{G}$ 
which themselves form simplicial complexes. Not all 
closed sets in the Whitney complex of a graph do have to be induced from a 
sub-graph of $G$. The boundary $S(x)$ of a maximal simplex for example is always
closed but it is only a {\bf skeleton complex} of a Whitney complex of a graph. 
\footnote{Taking subgraphs as closed sets is motivated by the Zariski topology,
where closed sets are algebraic subsets $H$ of an algebraic variety $G$.}
To see the correspondence between open and closed sets: note that any simplicial 
complex $K$  in $\mathcal{G}$ has as a complement the union of 
all stars $U(x)$ with $V(x) \cap V(K) = \emptyset$ and this union is an open set. 
For a graph $G$, a closed set $K$  contains all the 
simplicial complexes of {\bf subgraphs} of $G$. 

\paragraph{}
In the case $G=K_3$, where we have 
\begin{tiny}
$\mathcal{G}= \{ \{1,2,3\}, \{1,2\},\{2,3\},\{1,3\},\{1\},\{2\},\{3\} \}$. 
\end{tiny}
The topology $\mathcal{O}=\{U_1, \cdots U_{19} \}$ has 19 elements: \\
\begin{tiny} \parbox{13cm}{$U_{1}=\emptyset$,
$U_{2}$ = $\{(1,2,3)\}$,
$U_{3}$ = $\{(1,2),(1,2,3)\}$,
$U_{4}$ = $\{(1,3),(1,2,3)\}$,
$U_{5}$ = $\{(2,3),(1,2,3)\}$,
$U_{}$  = $\{(1,2),(1,3),(1,2,3)\}$, 
$U_{7}$ = $\{(1,2),(2,3),(1,2,3)\}$, 
$U_{8}$ = $\{(1,3),(2,3),(1,2,3)\}$, 
$U_{9}$ = $\{(1),(1,2),(1,3),(1,2,3)\}$,
$U_{10}$ = $\{(2),(1,2),(2,3),(1,2,3)\}$, 
$U_{11}$ = $\{(3),(1,3),(2,3),(1,2,3)\}$, 
$U_{12}$ = $\{(1,2),(1,3),(2,3),(1,2,3)\}$, 
$U_{13}$ = $\{(1),(1,2),(1,3),(2,3),(1,2,3)\}$,
$U_{14}$ = $\{(2),(1,2),(1,3),(2,3),(1,2,3)\}$,
$U_{15}$ = $\{(3),(1,2),(1,3),(2,3),(1,2,3)\}$, 
$U_{16}$ = $\{(1),(2),(1,2),(1,3),(2,3),(1,2,3)\}$, 
$U_{17}$ = $\{(1),(3),(1,2),(1,3),(2,3),(1,2,3)\}$, 
$U_{18}$ = $\{(2),(3),(1,2),(1,3),(2,3),(1,2,3)\}$, 
$\mathcal{G}=U_{19}=\{(1),(2),(3),(1,2),(1,3),(2,3),(1,2,3)\}$.} \end{tiny} \\
The number of closed sets $\{ U^c, U \in \mathcal{O} \}$ is of course the same than the 
number of open sets. The sub-graphs: there is the empty graph of dimension $-1$, 
7 graphs of dimension $0$, 3 graphs of dimension 1 with one edge, 3 graphs withe one edge and one vertex, 
3 graphs of dimension 1 with 2 edges and then
the complete graph. This gives $1+7+3+3+3+1=18$ sub-graphs. The $1$-dimensional skeleton
complex $C_3$ of $K_3$ is a closed set in the topology which does no correspond to a 
sub-graph. Its complement is the open set $U( (1,2,3) ) = \{ (1,2,3) \}$. 

\paragraph{}
Let $A$ be an arbitrary set of sets in $\mathcal{G}$. The {\bf closure} of $A$ is 
the smallest closed set (simplicial complex) which contains $A$. \footnote{This corresponds
to a classical notion in the theory of simplicial complexes when they are considered
as geometric realizations in Euclidean space. The closure contains all boundary simplices,
meaning to look at all subsets of $x$ and so also in the continuum it means
to look at the smallest simplicial complex which contains $K$.} 
For example, if $A = \{x\}$ consists of a single simplex, then its closure
$\overline{A}$ is the simplicial complex $\{ y \subset x \}$ generated by $A$. 
If $\mathcal{G}$ is the Whitney complex of a graph, the closure of a set of simplices
is often (but not always) the subgraph of the Whitney complex of the smallest sub-graph of $G$ 
which contains all simplices $x$. The example of the closed set 
$C_3 = \{ (1),(2),(3),(1,2),(2,3),(3,1) \}$ in $K_3$ which is the skeleton complex of the
Whitney complex of $K_3$ shows that not all closed subsets in the topology of $G$ 
correspond to sub-graphs of $G$. \footnote{One could define an other topology, where the closed sets
are {\bf simplicial complexes of sub-graphs} but this is a rougher topology and closer to the 
{\bf Zariski topology}.}

\paragraph{}
If $X$ is a topological space and $Y$ a subset, then $Y$ is called {\bf locally closed}
if it is an intersection of an open set $A$ and closed subset $K$. We think then 
of the induced topology of $X$ on $K$ by taking all intersections $U \cap K$ as 
open sets in $K$. By definition, every $U \cap K$ is an intersection of a closed and an
open set.  Locally closed sets are not 
necessarily closed in the topology $X$ but part of the Borel $\sigma$-algebra of the topology. 
Locally closed sets are sets which are open in some closed subset $K$ with induced topology. 
In our context, where a closed set is a simplicial complex, a locally closed set
is a set which is an open set in that simplicial complex. It does not need to
be open in the original topology $\mathcal{O}$. 
Take for example the set which consists of a single simplex $A=\{ x \}$ which is not 
maximal, nor zero-dimensional. An example is if $x$ is a boundary edge in a Wheel graph $G$. This set is 
neither open nor closed. But the set is locally closed because we can write it as an open 
set in the boundary sub complex $K$ which is a circular graph complex $K$. The set $A$ is an open set 
in $K$ but not open in $G$. The complement of $A$ in $K$ is closed in $K$ as well as closed in $G$. 
Similarly, look at the closure of $A$ which is the complete complex with $2$ elements. Its
complement in $K$ is open in $K$ but not open in $G$. 

\paragraph{}
Let us look at some examples: in a finite abstract simplicial complex, every single simplex 
$\{x\}$ is locally closed because it is the intersection of the open set 
$U(x)$ and the closed set $W^-(x)=\overline{ \{x\} \} }$. 
In the simplicial complex of $K_3$, the set $\{ \{1\},\{1,2\},\{1,2,3 \} \}$ is not
locally closed. From the $128$ possible subsets of $\mathcal{G}$, there are $64$ which are not 
locally closed and 64 which are. In a one-dimensional simplicial complex, all subsets are
locally closed. In $K_4$, where we have $2^{15} = 32768$ possible set of 
subsets in $\mathcal{G}$, where only $167$ sets are open and $167$ are closed, 
there are $3605$ locally closed sets.

\satz{The finite topology on a complex is of Zariski type:
sub simplicial complexes are the closed sets. 
Sub-graphs of a graph define closed sets of a graph but not all closed sets come from subgraphs. 
On complexes coming from graphs, we could get a slightly rougher topology by declaring 
sets to be closed if their simplicial complexes are Whitney complexes 
coming from closed sub-graphs. Finally, we have looked at locally closed sets, 
sets which are intersections of open and closed sets. }

\section{Compact} 

\paragraph{}
Traditionally, a topological space is declared to be {\bf compact} if every open cover has a 
finite sub-cover. Compactness in this strict sense is not a very useful in finite topological spaces:
using the definition, {\bf every set} is compact (whether it is open, closed and even if it is
neither). It makes therefore much sense in a finite topological space to to identify compact sets with
closed sets. If topologies are considered for non-finite graphs, 
a set then would be compact in a complex or graph, if it is a closed and finite set 
defined by a finite abstract simplicial complex. 
Again, it would make sense to consider {\bf being compact} as a synonym for being finite
{\bf and} closed. In an infinite graph $G$, many (but not all) compact sets are represented
by finite subgraphs of $G$. Their Whitney complexes are then simplicial complexes and also closed.
On a general graph $(V,E)$ with no restrictions on $V,E$, one 
could also look at the slightly rougher topology in which finite subgraphs define closed sets. 
There are less closed sets then and therefore also less open sets. 
In zero dimensions, we would get the {\bf co-finite topology},
which in the finite case agrees with the {\bf discrete topology} 
and where {\bf compactness} is equivalent to being {\bf finite}. 

\paragraph{}
We have just seen that in general, like if one looks at topological spaces which are 
not so commonly used - finite topological spaces are examples - one has 
to be a bit more careful when using notions which involve compactness. The property that every 
open cover has a finite subcover is really not a good notion for compactness in the 
case of finite topological spaces. For example, 
a map is called {\bf proper} if the inverse image of a compact set is 
compact. In the case of finite topologies with the standard definition of compactness, 
{\bf every map} (even a not continuous one)
would be proper, simply because the inverse image of any finite set is a 
finite set, which by the definition of having a finite
cover, would be declared to be compact. So, the notion of ``proper" does not really 
say much. If one requires a compact set to be closed too, then on finite topologies, 
the notion of {\bf proper} is the same than continuous because there, every 
closed set is compact and the fact that all inverse images of closed sets 
is closed is equivalent to continuity. 

\paragraph{}
Looking art finite topological spaces could be named ``radically elementary topology"
similarly as radically elementary probability theory covers a lot of traditional 
probability theory \cite{Nelson}. 
Finite structures are not that limiting, especially if one considers them in a 
{\bf non-standard frame work}. In internal set theory IST for example, 
\cite{Nelson77,NelsonSimplicity} compact sets can be
modeled by finite sets. Compact simplicial complexes $X$ therefore can be modeled by 
finite abstract simplicial complexes. Of course, the number of elements is non-standard
if the set $X$ is infinite. If one looks at a continuous map from a compact topological 
space to itself, then in general there are infinitely many fixed points. In a situation 
like in the context of the {\bf Lefschetz fixed point theorem}, one traditionally 
assumes that there are finitely many fixed points. 
In the non-standard frame work, one would just assume that
the number of fixed points is {\bf standard} (an axiomatically defined term). 
This will then assure that also the sum of the Lefschetz numbers is standard. 

\satz{In a finite topological space, every set is compact so that the classical notion is not
very useful. Every map would be proper for example. 
It makes sense therefore to consider all closed finite sets as compact instead. 
If one would consider infinite complexes, every finite 
sub-simplicial complex would be considered compact. Especially, every finite
subgraph graph defines a closed compact set. In non-standard analysis frame-works,
"compact topological spaces" can be treated like "finite topological spaces".
}

\section{Connectivity} 

\paragraph{}
A graph $G$ is called {\bf path connected} if for every two vertices $a,b$,
there is a {\bf path} $e_1=(a,v_1),\dots,e_n=(v_{n-1},b)$ (a finite collection of edges) 
connecting $a=v_0$ with $b=v_n$. A simplicial complex $\mathcal{G}$ is {\bf path connected}
if the graph $G_1$ defined by the complex is path connected. 
A graph is path connected if and only its topology $\mathcal{O}(G)$ is connected: 
here is the proof: if $\mathcal{G} = U \cup V$ is the disjoint union of two open sets $U,V$,
then there can not be simplex which contains $x \in U$ and $y \in V$. This means that 
there is no path connecting $x$ with $y$ in the graph $G_1$ and $\mathcal{G}$ is not path 
connected. On the other hand, if $K_1,K_2$ are not path connected components in $G_1$, then
their smallest open neighborhoods $U_1,U_2$ are disjoint and $\mathcal{G}$ is not connected. 

\paragraph{}
In graph theory, connectivity is also called {\bf $0$-connected}. 
A graph is called {\bf $1$-connected}, if it is connected but there exists a vertex $v$ 
which when removed, renders $G \setminus v$ disconnected. The graph $K_2$ is $1$-connected
but its Barycentric refinement, the path graph $P_3$ is not $1$-connected. 
A simplicial complex can be defined to be {\bf $1$-connected}, if is graph $G_1$ is $1$-connected. 
If $f: G \to H$ is a continuous surjective map between graphs and 
$G$ is $1$-connected, then also $H$ is $1$-connected: to prove this, note first 
that $H$ must be connected because the continuous image of a connected topological space is 
always connected. If $H$ was not 1-connected, all unit spheres $S(x)$ in $H$ would be connected.
So, also $f^{-1}(S(x))$ is connected.  
This does not generalize to higher connectivity in graph theory: 
it is possible that the continuous image of a $2$-connected graph can be
only $1$-connected. As an example, take the {\bf kite graph} $G$ in which one edge has
been removed from $K_4$. This graph is $2$-connected, but a
continuous map from $G$ to $K_2$ has an image that is only $1$-connected
but not 2-connected. 

\satz{The topology of a complex is connected if and only if the complex is
classically connected in the sense that the graph of the complex is path 
connected. The topology of a graph on its simplicial complex is connected if and only 
if it is path connected. A complex is connected if and only if its graph
is connected. A continuous image of a connected complex is a connected complex.
If a graph is $1$ connected, then a continuous image is still $1$ connected.}

\section{Separation}

\paragraph{}
The topology is called a {\bf Kolmogorov space} or $T_0$,
if for any pair of distinct points $x,y$, there is at least one point who has a neighborhood 
not containing the other. The topology $\mathcal{O}$ of a simplicial complex $\mathcal{G}$ 
is always a Kolmogorov space: this is clear in the case $x \subset y$ or $y \subset x$. 
In the case $x \subset y$, then $U(y)$ does not contain $x$. 
If $x \cap y$ is not empty, then $U(x)$ does not contain $y$ and $U(y)$ 
does not contain $x$. 

\paragraph{}
However, the topology $\mathcal{O}$ is not {\bf Fr\'echet} = $T_1$ if the dimension is positive. 
Fr\'echet means means that for every two points $x,y$, there exists a neighborhood
of $x$ which does not contain $y$ and a neighborhood of $y$ that does not contain $x$.
As an example, take $x \neq y$ with $x \subset y$. Now, only $y$ can be separated from $x$
but $x$ can not be separated from $y$. Every neighborhood of $x$ contains $y$. If the 
complex has positive dimension, the it is never Fr\'echet. Of course, all zero-dimensional 
complexes are Fr\'echet because the topology is the discrete topology. 

\paragraph{}
The topology of a simplicial complex $G$ is also not {\bf Hausdorff} or $T_2$
if the complex has positive dimension. In particular, the topology of a graph $G$ 
is not Hausdorff if $G$ has positive dimension. 
Two vertices $x,y$ in a graph which are connected by an edge can not be separated by open sets
in the complex: as x must contain $U(x)$ and $y$ must contain $U(y)$ and these two open sets 
have an intersection $U(x) \cap U(y)$ which contains the edge $\{e\}$. 

\paragraph{}
The topology is also not {\bf normal} = $T_4$: two closed sets can in general not be separated by 
open sets. Examples are the closures $\overline{x},\overline{y}$ of two simplices
$x,y \in \mathcal{G}$ that have a non-empty intersection. 
They are closed sets but every neighborhood of one intersects with any neighborhood 
of the other. One could think that the non-Hausdorff property is a handicap. However, 
it can be a blessing in the context of {\bf connection calculus}, when we consider
higher order {\bf characteristics} generalizing the Euler characteristic 
$\chi(G) = \sum_{x \in \mathcal{G}} \omega(x)$ with $\omega(x) = (-1)^{\rm dim}(x)$. While
Euler characteristic satisfies for every subset $A,B \subset \mathcal{G}$ the {\bf valuation 
property} $\chi(A \cup B) = \chi(A) + \chi(B) - \chi(A \cap B)$, this property does no more 
hold in the case of {\bf Wu characteristic }
$\omega(G) = \sum_{x,y, x \cap y \in \mathcal{G}} \omega(x) \omega(y)$. But the valuation 
property will hold for open sets.

\paragraph{}
A topology is called {\bf Alexandroff} if every point $x$ has a smallest non-empty open 
neighborhood $U(x)$. In an Alexandroff topology, space has {\bf smallest atoms}. 
Equivalently, an Alexandroff topological space has the property that arbitrary intersections
of open sets are open. Most topological spaces we are familiar with are not Alexandroff. 
If a metric space is Alexandroff, it must be a discrete topological space as then, 
every single point needs to be open. In general, any discrete topology and the 
{\bf indiscrete topology} $\mathcal{O}=\{ \emptyset,X \}$  is Alexandroff. 
Any {\bf finite topology} must be Alexandroff because the intersection of all open sets containing
$x$ is an open set $U(x)$, containing the point. It is the {\bf star}. 
Calling a finite topology a ``finite Alexandroff topology" would be a pleonasm. Still, it is good
to use the name as Alexandroff was one of the firs who seriously considered finite topological
spaces.

\paragraph{}
To compare, note that the geometric realization $|G|$ (in some Euclidean space) of a
finite abstract simplicial complex $\mathcal{G}$ (or finite simple graph 
with the Whitney complex) is always a Hausdorff, because $|G|$ is a closed subset of a
Hausdorff topological space. We will see that looking at the
topological realization of a complex loses some information like the 
topological nature of spheres in the space. 
Topological manifolds can in general not be described by one simplicial 
complex (or equivalence class of Barycentric refinements) alone. This
is not a surprise. Most topological spaces we look at, even compact ones
like Cantor type sets can not be described by one single {\bf finite 
standard topological space}. 
\footnote{Nonstandard analysis teaches us however that
we can describe it by a non-standard finite topological space. Standard 
finite topological spaces are then the spaces we look at here. The axiom
system however does not allow us to define the intersection of all standard
open sets so that compact topological spaces do not have atoms $U(x)$.}

\satz{The topology of a complex or graph is Kolmogorov ($T_0$) not Fr\'echet (not $T_1$), 
not Hausdorff (not $T_2$) and not normal (not $T_4$), but (like all finite topological spaces) 
is Alexandroff. The non-Hausdorff property is in sharp contrast with the topology given 
by geometric realizations which are Hausdorff. }

\section{Dimension} 

\paragraph{}
The {\bf maximal dimension} of a simplicial complex is defined as
${\rm max}_{x \in \mathcal{G}} {\rm dim}(x)$.
The {\bf maximal dimension} of a graph $G$ is the maximal dimension of its Whitney complex $\mathcal{G}$. 
We have already seen that a continuous map $f: \mathcal{G}$ to $\mathcal{H}$ between simplicial
complexes has the property that ${\rm dim}(f(x)) \leq {\rm dim}(x)$. The maximal dimension of a 
continuous image of a complex $\mathcal{G}$ is therefore smaller or equal than the dimension of $\mathcal{G}$. 
An {\bf open cover} $\{ U_j \}$ of $G$ is a set of open sets such that 
$\bigcup_{j} U_j = \mathcal{G}$. A cover defines a \v{C}ech {\bf nerve graph}, in which the sets $U_j$
are the vertices and where two are connected if they simultaneously intersect in a non-empty set. 
The maximal dimension of this graph is called the {\bf dimension of the cover}.
The minimum over all dimensions of covers of $G$ is called the {\bf topological dimension} of $G$. 

\paragraph{}
The following fact is an other reason why the topology associated to a simplicial complex or 
graph is the right one. 
The topological dimension and the maximal dimension $d$ are the same for every complex
and every graph. The reason is that we can cover the space with sets $U(\{v\})$ with 
$v \in V = \bigcup_{x \in \mathcal{G}} x$. This cover can 
not be refined further because removing one would keep some vertex $v$ uncovered. 
The dimension of this cover is equal to the maximal dimension of $G$ because if $x$ is a simplex
of dimension $d$, then all the open sets $\{ U(v) \}_{v \in \mathcal{G}}$ intersect. 
Therefore, the topological dimension of $G$ is at at most $d$. Every cover of $G$ must 
contain all sets $U(\{v\})$ with $v \in V$ because otherwise $\{v\} \in \mathcal{G}$ 
would not be covered. So, the dimension is also at least $d$. 

\paragraph{}
The \v{C}ech graph of an open cover $\mathcal{U}$ of a graph or simplicial 
complex is defined as the graph in which the open sets $\mathcal{U}$ are the 
vertices and where two vertices are connected if they have a
non-empty intersection. In general, the \v{C}ech graph of the cover $\{ U( \{ v\} ), v \in V \}$ is the 
graph $G$ itself and the \v{C}ech cover of the cover $\{ U(x), x \in \mathcal{G} \}$ 
is the Barycentric refinement $G_1$. One usually looks at covers for which every of the open 
sets are contractible (every star $U(x)$ is considered contractible because its closure
$B(x)$ is contractible.)

\satz{The topological dimension of a topology on a complex as defined by Lebesgue, 
agrees with the maximal dimension of the complex. The \v{C}ech graph of the
base cover of a complex is the graph $G_1$. If $\mathcal{G}$ is the Whitney 
complex of a graph $G$, then the \v{C}ech graph of the cover $\{ U(v) \}_{v \in V}$
(where $V$ is the set of zero dimensional simplices) is the graph $G$ itself. 
}

\section{Product}

\paragraph{}
Every data structure, whether we deal with graphs, with simplicial complexes or with 
topological spaces has notions of products. 
The {\bf Shannon product} $G*H$ \cite{Shannon1956} of two graphs $G,H$ is the graph for which 
$V(G*H) = V(G) \times V(H)$ is the {\bf Cartesian product} of sets and where 
$E(G*H)=\{ ((a,b),(c,d)), (a,c) \in E(G) \cup V(G) \; {\rm or} \; (b,d) \in E(H) \cup V(H)  \}$,
meaning that two points are connected if both projections have the property that 
they project onto a vertex or edge.  The Shannon product does not go over naturally to 
complexes. We can of course look at the complex of the Shannon product. 
What is nice about the Shannon product is that it allows to see graphs as a {\bf ring}.
We have explored this a bit in 
\cite{ArithmeticGraphs,RemarksArithmeticGraphs,ComplexesGraphsProductsShannonCapacity}.

\paragraph{}
The {\bf box product topology} 
\footnote{For finitely many products, the box product topology agrees with the 
product topology.}
on $G*H$ is the finest topology on the simplicial complex of $G*H$ such that
both the projections on $G$ and $H$ are both continuous. 
The graph topology of the product $G*H$ is in general much finer than the product topology. 
One can see this already for $G=H=K_2$, where $G*H=K_4$. The topologies of $G$ and $H$
have only $5$ elements, while the topology of $G*H$ has $167$ elements. 

\paragraph{}
The graph topology on $U*V$ is in general finer than the topology generated by the 
``cubes" $U*V$ where $U,V$ are basis elements of the factors. 
Some simplices in $G*H$ are of the form $x*y$ which is a $(k+1)*(l+1)-1$-simplex if 
$x$ was a $k$ simplex and $l$ was a $l$ simplex. 
But not every simplex in $G*H$ is of the form $x*y$. For example, 
for $G=H=K_2$, only the $0,1$ and $3$-dimensional simplices in $G*Y$ are products. 
The 2-simplices (the triangles) in $G*H$ are not products because $2=(k+1)*(l+1)-1$ implies 
either $k+1=3$ or $l+1=3$ but there are no 2-simplices in neither $G$ nor $H$. 
One can see this also from the fact that there are $2^2-1=3$ simplices in $G$ and 
$H$ and $2^4-1=15$ simplices in $G*H$. Only 9 of them are of the form $x*y$ with 
$x$ a simplex in $G$ and $y$ a simplex in $H$. 

\paragraph{}
The Shannon product does not preserve manifolds. 
The {\bf Stanley-Reisner product} of two simpicial complexes $\mathcal{G}$ and $\mathcal{H}$
is defined as the Whitney complex of the graph in which the Cartesian product
$\mathcal{G} \times \times \mathcal{H}$ are the vertices and where two different 
vertices $(x,y),(u,v)$ are connected by an edge if either $x \subset u, y \subset v$ or 
$u \subset x, v \subset y$. The Stanley-Reisner product of a $p$-manifold with a $q$-manifold
is a $(p+q)$-manifold. As the graph defining $\mathcal{G} \times \mathcal{H}$ is homotopic
to the Shannon product, it inherits properties of the former like the K\"unneth formula
or the compatibility with the Euler characteristic. We will write more about the compatibility
with higher characteristics elsewhere and especially show that they are topological invariants. 

\paragraph{}
One can now ask whether there is not a natural ring structure on simplicial complexes which 
corresponds to the Shannon ring, or whether there is even a ring structure which preserves
manifolds. One problem is that if we take the Cartesian product of two simplicial complexes,
we don't have a closed set. We would have to close it but it would not have the properties we like.
Much more elegant is to expand the class from simplicial complexes to {\bf delta-sets}.
This structure is more general than {\bf simplicial sets}, a popular construct which has more
structure than $\delta$-sets. Every simplicial set of course is also a delta set by just forgetting
the degeneracy maps $s_i$ maps and only keep the face maps $d_i$. The disadvantage of working
with $\delta$ sets is that we have to carry around not only sets of sets but also keep track 
of the  maps $d_i$. And the entire elegance of having a simple topology etc is gone. 
$\delta$ sets also are useful when describing {\bf quivers}, which generalize finite simple graphs. 
There is a Whitney functor from quivers to $\delta$ sets generalizing the 
functor from graphs to simplicial complexes. 
\footnote{In order to talk about funtors one needs to adapt the morphisms on graphs and
simplicial complexes and continuous maps are the most natural common denominator as
both simplicial maps as well as graph homomorphisms are continuous maps.}

\satz{The projections from the Shannon product $G*H$ of two graphs to one of its factors
is a continuous map. The graph topology of the product graph is in general much finer than 
the product topology in general. The Shannon product does not preserve topological quantities
like higher characteristic or manifold properties. The Stanley-Reisner product however does. 
The Stanley-Reisner product is more compatible with topology but does not provide a ring
structure as associativity fails. The Shannon product on the other hand defines a ring and so
an arithmetic. 
}

\section{Join}

\paragraph{}
The {\bf join} $G \oplus H$ of two graphs $G,H$ has as vertex set $V(G) \cup V(H)$ and
as the edge set $E(G) \cup E(H) \cup \{ (a,b), a \in V(G), b \in V(H) \}$. The join 
operation in graphs theory was first defined by Zykov  
\cite{Zykov} and does exactly what the join does for geometric realizations of the
Whitney complex. The join with a $0$-sphere is a {\bf suspension}. 
In general, the join of two spheres is a sphere. The proof follows from the sphere
formula $S_{G \oplus H}(x) = S_G(x) \oplus H$ and $S_{G \oplus H}(y) = G \oplus S_H(y)$. 
so that by induction if each $G,H,S_G(x),S_H(y)$ are spheres, then the unit sphere
of any point in $G \oplus H$ is a sphere. 
The join of a graph with the $0$-sphere $S_0$ is called the {\bf suspension}. 
Since $S_0 \oplus S_0 = C_4$ is a cyclic graph and so a discrete sphere, the join of a graph 
with a cyclic graph is a {\bf double suspension}. 

\paragraph{}
The {\bf join of two simplicial complexes} $\mathcal{G},\mathcal{H}$ can be defined also
without referring to their graphs. The simplices in $\mathcal{G} \oplus \mathcal{H}$
are the union of $\mathcal{G}, \mathcal{H}$ and $\mathcal{G} \oplus \mathcal{H}$
where $x \oplus y = x \cup y$ is a $k+l+1$-dimensional simplex obtained by taking
the disjoint union of the two sets $x,y$. If $x$ had $k+1$ elements and $y$ had $l+1$
elements, then $x \oplus y$ has $k+1+l+1=k+1+2 = (k+l+1) + 1$ elements so that 
$k+l+1$ is the dimension of $x \oplus y$. So, for example, if 
$\mathcal{G} = \{ a,b \}$ is the zero sphere and 
$\mathcal{H} = \{ 1,2,3,4,(12),(23),(34),(41) \}$
is a discrete circle, then $\mathcal{G} \oplus \mathcal{H} = \} 
=  \{ a,b, 1,2,3,4,(12),(23),(34),(41), (a1),(a2),(a3),(a4),(a12),(a23),(a34),(a41)$ \\
   $(b1),(b2),(b3),(b4),(b12),(b23),(b34),(b41) \}$ is the suspension of a circle and the 
{\bf octahedron complex}. 

\paragraph{}
The join operation is {\bf dual} to the {\bf disjoint union} $G + H$ as addition: 
if $G'$ denotes the graph complement of $G$ in which edges and non-edges are switched, then 
$(G \oplus H)' = G' + H'$, where $G+H$ is the 
{\bf disjoint union} of the graphs. If $G_1$ is the Barycentric refinement of $G$, 
then the unit sphere $S(x)$ of a point $x \in \mathcal{G}=V(G_1)$ is of the 
form $S^+(x) \oplus S^-(x)$, where $S^+(x) = \{ y, x \subset y\}$ is the 
{\bf unstable sphere} and $S-(x) = \{ y, y \subset x \}$ is the {\bf stable sphere}. 
In a discrete manifold, where every unit sphere is a sphere, 
both the stable and unstable spheres are spheres. 

\paragraph{}
If $G,H$ are graphs and if $x$ is a simplex in $G$ and $y$ is a simplex in $H$, 
then $x \oplus y$ is a simplex in $G \oplus H$. So, the Whitney complex of the Zykov join 
$\mathcal{G}(G \oplus H) = \mathcal{G} \oplus \mathcal{H}$ is the join of
the simplicial complexes which is the 
$\mathcal{G} \cup \mathcal{H} \cup \{ x \oplus y, x \in \mathcal{G}, y \in \mathcal{H} \}$. 
The topology of the join by definition has as a basis the sets 
$U_{G \oplus H}(z)$, where $z=x \oplus y$ is a simplex in $G \oplus H$. 
This is $U_{G \oplus H}(x \oplus y) = U_{G \oplus H}(x) \oplus H \cup G \oplus U_{G \oplus H}(y)$. 
The total set of stars  $U(x \oplus y)$ is a basis of $G \oplus H$ and generates
the topology of $G \oplus H$. 

\satz{The embedding of $G$ in $G \oplus H$ with the induced topology is a classical 
homeomorphism onto the image, as it is a bijection onto the image. 
The join of the topological base in $G$ and $H$ defines a base for the join $G \oplus H$. 
The smallest atomic open sets  $U(x \oplus y)$ in $G \oplus H$ is a basis: the base
of the join contains open sets $U(x) \oplus H$ as well as open sets $G \oplus U(y)$. 
}

\section{Subgraph}

\paragraph{}
We have seen that a subcomplex $\mathcal{H}$ of a complex $\mathcal{G}$ is a closed set. 
The topology induced from $\mathcal{G}$ on a subcomplex $\mathcal{H}$ 
is the topology of $\mathcal{H}$ itself. A general subset $\mathcal{H}$
of $\mathcal{G}$ can still be given a topology by taking as open sets $U \cap \mathcal{H}$ with 
$U$ in the topology of $\mathcal{G}$. This topological space agrees with the 
closure of $\mathcal{H}$. Similar standard consequences hold for graphs. 
Any subgraph $H$ of $G$ has a topology which agrees with the induced 
topology from $G$. A subset $A$ of $\mathcal{G}(H)$ is open if and only if it is of the form 
$\mathcal{O} \cap \mathcal{G}(H)$. The simplicial complex of a subgraph $A$ of $G$
is a closed set. As in general when we take the relative topology on a closed subset $K$ of a
topological space $X$, the {\bf relative topology} has the sets $U \cap K$ as open sets, where
$U$ ranges over the open sets in $G$. All the axioms for a topology are satisfied. The relative
topology is the finest topology on $K$ which has the property that the inclusion $i:K \to G$ is
continuous. 

\paragraph{}
To reformulate this, the relative topology on a subgraph $H$ of $G$ agrees with the 
graph topology of $H$. We only have to look at a basis to see this. If $x$ is a simplex in $H$,
then $U_H(x)$ is the intersection $U_G(x) \cap H$.  So, the basis for the topology on $H$
is the same than the restriction of the basis of the topology on $G$ to $H$. This proves the 
statement. The subgraph $H$ can be generalized. As for any topological space on some set $X$
we can restrict the topology on any subset $Y$ of $X$. So, we can build a topology on any 
subset of the simplicial complex $\mathcal{G}$ of a graph $G$. 

\paragraph{}
If a subset $W \subset V$ of the vertex set of a finite simple graph $G=(V,E)$ is given, 
one can look at the subgraph $H$ {\bf generated} by $W$. This means to take the largest
subgraph of $G$ which contains the vertex set $W$. One could look at the closure of a
subgraph $H$ as the subgraph generated by the vertex set of $W$. This is in general 
a much larger graph. For a {\bf Hamiltonian subgraph} of $G$ (a graph which passes
through all vertices) for example this ``closure" would be the graph itself. The topology
defined by the simplicial complex however make a subgraph naturally closed already
as the simplicial complex $\mathcal{H}$ generated by the simplices in $H$ is a sub-simplicial complex
of the complex $\mathcal{G}$ of $G$ and so closed. 

\satz{Sub simplicial complexes correspond to closed sets in the topology. 
Subgraphs of a graph define a subclass of closed sets in the topology on the simplicial 
complex defined by the graph. The relative topology on a subgraph of a graph is the 
graph topology of $H$ itself and does not use the topology of the host graph $G$. 
The relative topology on a subcomplex of a complex is the topology of the subcomplex
itself which is the same without looking at the ambient space $\mathcal{G}$. 
}

\section{Quotient}

\paragraph{}
If $G$ is a graph and $\sim$ is an equivalence relation on vertices honoring the
edges, then the set of equivalence classes $H=G/\sim$ can carry a topology. First 
of all, the equivalence relation induces an equivalence relation on complete subgraphs
and $x/\sim$ is the complete graph on the set of equivalence classes of $V(x)/\sim$.
Define $\mathcal{G}(H)$ as the set $\{ x/\sim, x \in \mathcal{G}(G) \}$. 
If $G=K_2$ with $\mathcal{G}=\{ \{1,2\},\{1\},\{2\} \}$ for example 
and $\sim$ identifies the two points $1,2$, we get $\mathcal{G}(H) = \{ \{1\} \}$. 
Assume $G$ is a {\bf cover} of $H$, meaning that there is a surjective graph homomorphism $f: G \to H$, 
then we can see $H$ as a quotient $G/\sim$, where $v \sim w$ if $f(v)=f(w)$. An example
is the cover $C_8 \to C_4$ with $f(v)=v \; {\rm mod} \;  4$ if $V(C_8)=\{0,1,2,3,4,5,6,7 \}$. 
An other example is the cover $S^2 \to P^2$ of a sufficiently large $2$-sphere for which 
the equivalence relation defines a manifold. Some 2-spheres are too small. For an octahedron
$O$ for example, a graph with $6$ elements, identifying opposite vertices produces not
no projective plane but $O/\sim = K_3$. 

\paragraph{}
In order that an equivalence relation on the sets of simplicial complex $\mathcal{G}$  produces a 
a quotient topology, one needs to make some assumptions. In a Barycentric refinement, things are easier.
Having the complex $\mathcal{G}$ too small can make things  weird. Lets look for example 
the cycle complex $\mathcal{C}_4 = \{ 1,2,3,4,(12),(23),(34),(41) \}$
and impose the equivalence relation where we identify the vertices $1$ and $3$. The quotient is no
more a simple graph, but a {\bf quiver} because multiple connections appear. We have now a graph 
with three vertices $1,3,4$ and double bond connections $(13)$ and double bond connections $(34)$. 
We can however look at the situation in the Barycentric refinement, where the identification 
becomes now a figure $8$ graph. This situation matters if we look at {\bf Riemann-Hurwitz} formulas
which relate the Euler characteristic of the quotient with the Euler characteristic of the complex
itself as well as using {\bf ramification points}. In a case of a covering and having a group $A$
of order $|A|$ acting on $\mathcal{G}$ so that $\mathcal{G}/A$ is again a complex, then 
The Riemann-Hurwitz formula tells $\chi(G/A) = \chi(G)/|A|$. For example, if $A=\mathbb{Z}_2$ acts
on a sphere $\mathcal{G}$ and $\mathcal{H}=\mathcal{G}/A$ is a projective space, 
then $\chi(\mathcal{H}) = \chi(\mathcal{G})/2$. But if the complex is too small like for the
Octahedron complex considered above, then $\mathcal{G}/A$ is not a complex any more. For the 
Barycentric refinement however it works and we get like that complexes representing a projective
plane. 

\satz{The topology on a quotient $\mathcal{H}=\mathcal{G}/\sim$ can be defined as usual in topology as the 
finest topology which makes the projection from $\mathcal{G}$ to the space of equivalence classes
continuous. If $\mathcal{H}$ is a simplicial complex, then its topology is the quotient topology.
If quotients come from covers, then we have completely analogue situations like in the continuum.
We can for example look at the quotient of an antipodal map on a Barycentric refined sphere and
get a finite topological space representing a projective space. 
}

\section{Manifold}

\paragraph{}
A graph $G$ is called a {\bf $d$-manifold} if every 
unit sphere $S(x)$ in $G$ is a $(d-1)$-sphere. 
A {\bf $d$-sphere} is a $d$-manifold such that for some vertex $v$, the graph $G-v$
without $v$ is contractible. A graph is {\bf contractible} if there exists 
$v$ such that $S(v)$ and $G-v$ are both contractible. These notions can be
defined also for complexes without referring to graphs. 

\paragraph{}
A complex $\mathcal{G}$ is a {\bf $d$-manifold} if every unit sphere 
$S(x)=B(x) \setminus U(x)$ is a $(d-1)$-sphere, where $B(x)=\overline{U(x)}$ is the unit ball
the closure of $U(x)$. The unit $S(x)$ is always closed and so carries a simplicial 
complex structure. A $d$-sphere is a $d$-manifold $\mathcal{G}$ such that 
$\mathcal{G} \setminus U(x)$ is contractible for some $x$. A complex is 
{\bf contractible} if there exists $x$ such that $S(x)$ and $\mathcal{G} \setminus U(x)$ 
are both contractible. One can extend contractibility to non-closed sets by defining
for example an open set to be contractible if its closure is contractible. Every star
$U(x)$ is contractible with this definition. 
\footnote{We identify {\bf collapsible} and {\bf contractible} and use 
{\bf homotopic to 1} if we mean that a complex can be deformed to $1=K_1$  by both 
expansions and contraction steps. While homotopic to $1$ is a computationally 
difficult equivalence relation, contractibility is easy to check.} 

\paragraph{}
All these inductive definitions are primed by the assumption that the 
empty graph $0$ or the empty complex $0$ is the $(-1)$-sphere and that the 
one point graph $1=K_1$ or $1$-point complex $1$ is contractible. 
If $G$ is a graph that is a $d$-manifold, then all its Barycentric 
refinements $G_n$ are $d$-manifolds too. If $\mathcal{G}$ is a complex which is a 
$d$-manifold, then all their Barycentric refinements are manifolds.
Examples of discrete manifolds in the sense just defined are {\bf combinatorial
triangulations} of a manifold. But not all triangulations are manifolds. The tetrahedron $K^4$
for example is contractible and so not a manifold. We could look at the $2$-skeleton complex
of $K_4$ however and get a sphere complex. The Barycentric refinement of $K_4$ is a $3$-ball
with $2$-dimensional boundary which corresponds to the Barycentric refinement of the 
$2$-skeleton complex. 

\paragraph{}
Let us add a remark coming from the continuum:
every PL manifold (a manifold equipped with a PL-structure) admits a combinatorial triangulation. 
The question of Poincar\'e from 1899, whether every smooth manifold 
admits a triangulation has been answered positively in the 1930ies: every smooth manifold has an 
essentially unique PL-Structure. (The converse is not true. There are PL-manifolds which can not
 be smoothed or admit different smooth structures.) The question shifted then to the topological situation.
See \cite{HistoryTopology} and especially \cite{Scholz} for history or \cite{Manolescu} for more
recent developments. 
Also PL structures do not exist in general on topological manifolds in dimensions 4 or larger. 
4-manifolds examples were given using the Kirby-Siebenmann invariant. 
Non-PL triangulations of manifolds were constructed using the Edwards-Cannon 
double suspension theorem. 

\paragraph{}
The fact that the Barycentric refined graph $G_1$ obtained from the Whitney complex $\mathcal{G}$ of a 
$d$-manifold graph $G$ is a $d$-manifold can be proven by induction with respect to dimension.
Indeed, $G$ is a manifold if and only if $G_1$ is a manifold. 
If $G$ was not a manifold, then some unit sphere $S(v)$ in $G$ would not be a sphere. 
But then also $S_{G_1}(v)$ which is the Barycentric refinement of $S(v)$ would not be a sphere. 
By induction assumption (unit spheres are one dimension smaller) this is a contradiction.

\paragraph{}
If $G$ and $H$ are homeomorphic and $G$ is a manifold and $H$ is not, then every $G_n$ is a manifold
and non of the $H_m$ are. Take a unit sphere $S(v)$ in $H$. Since its inverse image is homeomorphic to $S(v)$
it is a sphere. By definition, it is the boundary of a ball. All unit spheres of interior points in this
ball are by definition spheres. Now, every vertex in $\mathcal{G_n}$ is in the 
interior of the inverse of a ball $B(x)$ in $\mathcal{H}$.
The following statement in the summary is not true if ``homeomorphic" would be replaced by 
``has a homeomorphic geometric realizations". The notion of homeomorphism proposed here
is probably is equivalent to PL-equivalent but we do not prove this because we don't deal with 
infinity here: 

\satz{If $\mathcal{G},\mathcal{H}$ are homeomorphic and 
$\mathcal{G}$ is a $d$-manifold then $\mathcal{H}$ is a $d$-manifold.}

\section{Contractible}

\paragraph{}
The concept of ``contractible" entered in a crucial way in the definition of ``sphere" and so in
the definition of manifold. The notion 
\footnote{Again: we avoid the term collapsible used often in the literature. } 
makes sense for general graphs and general simplicial complexes. It is different from 
{\bf homotopic to 1}, where one can do both homotopy extensions and contractions. 
The {\bf dunce hat graph} is a concrete example of a finite simple graph which is 
homotopic to a point which is not contractible. 
\footnote{The Dunce hat can be realized as a graph $G$ with 17 vertices, 
52 edges and 36 triangles. Its unit spheres are all either 1-spheres or
homeomorphic to figure eight graphs (wedge sums of two 1-spheres). 
There are homotopy expansions which make it contractible.}
The graph $G$ is contractible if and only if $G_1$ is contractible. 
A complex $\mathcal{G}$ is contractible if and only if $\mathcal{G}_1$ is contractible. 
Contractibility can be extended to non-closed sets by assuming the closure to be contractible. 

\paragraph{}
The continuous image of a contractible graph does not need to be contractible: an example 
is the map from the linear graph $L_{n+1}$ to the graph $C_n$ 
mapping the initial and end point to the same point.  While $L_{n+1}$ is contractible, the graph $C_n$ is not. 
It is also not true that if $H=f(G)$ is contractible then $G$ is contractible. 
Let $f$ map a $0$-sphere $G=\{ V=\{a,b\},\{\} \}$ to $K_1 = \{ \{a\},\{ \} \}$. This is continuous
because both graphs have the discrete topology but the $0$-sphere is not contractible. 

\satz{As in the continuum, homotopy transformations are not continuous in general. But 
contractions = homotopy reductions of graphs $f: G \to G-v$ or complexes 
$\mathcal{G} \to \mathcal{G} \setminus U(x)$ are continuous. In the 
Unlike the unit spheres $S(x)$, the balls $B(x)$ are always contractible.}

\section{Boundary}

\paragraph{}
A {\bf $d$-manifold with boundary} is a graph or complex with the property that every unit sphere is 
either a $(d-1)$-sphere or is a $(d-1)$ ball. The same definition applies for 
simplicial complexes. If $S(x)$ is a $(d-1)$ sphere, we have an {\bf interior point} $x$, 
if $S(x)$ is a $(d-1)$-ball, $x$ is a {\bf boundary point}.
For a manifold with boundary, the {\bf boundary} is a $(d-1)$ manifold without boundary. 
An example of a manifold with boundary is
a $d$-{\bf ball} which by definition is a $d$-sphere with a point removed. An other example is a 
complete graph $K_{d+1}$, where all points are boundary points. 
We can include $K_{d+1}$ into the class of 
manifolds with boundary just in order to have Barycentric invariance. We want a graph to be 
homeomorphic to its Barycentric refinement and in general to have a graph $G$ homeomorphic 
to $H$ if there exist continuous maps $G_n \to H_m \to G$. We can say:

\satz{If $\mathcal{G}$ is a manifold with boundary and $\mathcal{G}$ is 
homeomorphic to $\mathcal{H}$, then $\mathcal{H}$ is a manifold with boundary.
The same holds for graphs. }

\section{Duality}

\paragraph{}
Continuity is not compatible with some duality notions in graph theory or the theory
of simplicial complexes. 
The operation of mapping a graph to its {\bf graph complement} is in general not continuous. 
Already the dimensions do not work. The graph complement of a cyclic graph $C_n$ is always 
homotopic to a sphere or then to a wedge sum of two spheres \cite{GraphComplements}. 
Only for special cases like $C_5$, where the graph complement
is the same graph, the complement operation can be made to be a homeomorphism. 

\paragraph{}
One can ask however whether the graph complement operation maps homeomorphic graphs to homeomorphic
graphs. But also here the answer is no: let $f:G \to H$ be a continuous map between finite graphs 
like for example $f(x)=x \; mod \; 4$ from $G=C_8$ to $H=C_4$. 
Does there exist a continuous map from $G^c$ to $H^c$? 
The graph complement of $C^8$ is homotopic to a 2-sphere. The graph complement of $C_4$ is a 
disconnected union of two graphs $K_2$. Since $H^c$ is disconnected and $G^c$ is connected and
a continuous map can not map a connected space to a disconnected space (as continuity preserves
the property of being connected), we can also not see the graph complement operation as a map
preserving homeomorphic graphs.

\paragraph{}
Let us reformulate this a bit differently. 
We have identified simplicial complexes or graphs which are Barycentric
refinements of each other. The notion of graph complement is not at all compatible with the
Barycentric refinement notion. The Barycentric refinement of $C_4$ is $C_8$. But there is no
topological similarity between $C_4^c$ which is homotopic to a $0$-sphere and $C_8^c$ which 
is homotopic to a $2$-sphere. 

\paragraph{}
For simplicial complexes, there is the {\bf Alexander duality} operation: 
if $\mathcal{G}$ is a complex and $V \bigcup_{x \in \mathcal{G}} x$ 
is the set of  0-dimensional simplices, then the {\bf Alexander dual} of $\mathcal{G}$ 
is the complex $\{ y \subset V, (V \setminus y) \notin \mathcal{G} \}$. 

\paragraph{}
one can ask whether there is a duality notion which corresponds to the 
graph complement. If $G$ is a graph, we can look at the simplicial complex
$\mathcal{G}$ of $G^c$ but that does have little to do with the simplicial 
complex $\mathcal{G}$. One can experiment with other notions like
if $\mathcal{G}$ is an arbitrary simplicial complex and $V=\bigcup_x x$ is the 
vertex set of $\mathcal{G}$. We can look at the complement $\mathcal{G}^c$ of $\mathcal{G}$
in the complete complex on $V$. This is a duality notion, but of course, $\mathcal{G}^c$ 
is almost never a simplicial complex and the closure $\overline{\mathcal{G}^c}$
of $\mathcal{G}$ is a complex but $G \to \overline{\mathcal{G}^c}$ is not a duality
notion. We see that graphs have their purpose especially with respect to arithmetic.

\satz{While interesting for other reasons, duality notions like 
graph complement or Alexander duality are not 
compatible with homeomorphisms. The topic of duality also shows that
having different data structures for finite geometries is useful. 
The graph complement duality for graphs for example works well with 
arithmetic. It provides an isomorphism between the Sabidussy ring 
with join and large product as operations,
and the Shannon ring with disjoint union and Shannon product as
operations. 
}

\section{Edge refinement}

\paragraph{}
{\bf Edge refinement} are topological transformations of graphs. They induce 
topological modifications of simplicial complexes. It is best described on graphs.
Given a graph $G$ and an edge $e=(a,b)$, we can refine the graph by 
adding a new vertex $v$, remove $e$ and connect $v$ to the 
intersection of $S(a)$ and $S(b)$. When applied to the cyclic graph 
$C_n$, it produces $C_{n+1}$. 
When applied more generally to a discrete $d$-manifold, we get a new
$d$-manifold. The reason is that the new unit sphere $S(v)$ is the join of
the $0$-sphere $\{ a,b \}$ and the $(d-2)$-sphere $S(a) \cap S(b)$ and 
so again a $(d-1)$-sphere. The unit spheres $S(a)$ and $S(b)$ are not changed.
The unit spheres of vertices $w$ with the edge $e=(a,b)$ in the unit sphere 
are themselves edge refined. Using induction with respect to dimension,
one has now verified that edge refinements preserve basic invariants 
like Euler characteristic.

\paragraph{}
In full generality, an edge refined graph $G_e$ is homeomorphic to $G$. 
First we can define a surjective continuous map from $G_e$ to $G$ induced
from the rule that every vertex goes itself and that the new vertex goes to 
$a$. In order to see that this map is continuous, check the properties. 
The inverse of every star in $G$ is homeomorphic to a star in $G_e$.
The star $U(x)$ of $x=\{a\}$ has as as an inverse image the union of
the star of $a$ and star of $v$. The star of the edge $(a,v)$ has the empty
set as an inverse image, The star of the edge $(v,b)$ has as an inverse image 
the star of $(a,b)$. Similarly, the star of any simplex containing $(a,v)$ 
has an empty inverse image while the star of any simplex containing 
$(v,b)$ has as the inverse image the star of the corresponding simplex 
containing $(a,b)$. 

\paragraph{}
To check the other direction, we have to construct a continuous
map from the Barycentric refinement $G_1$ of $G$ to $G_e$. We can 
take the canonical homeomorphism projection map from $G_1$ to $G$ and 
modify it so that it becomes a map from $G_1$ to $G_e$. 

\paragraph{}
We should mention that also the Dehn-Sommerville property is 
preserved by edge refinement and Barycentric refinements.  
Dehn Sommerville spaces generalize spheres and like spheres produce a monoid
under the join operation. They therefore can be used to generate spaces
more general than manifolds but still have many properties of manifolds like
that odd-dimensional manifolds have zero Euler characteristic. 
{\bf Dehn-Sommerville $d$-spaces} $\mathcal{X}_d$ \cite{DehnSommerville} are 
inductively defined. They must have the property that $\chi(G)=1+(-1)^d$ 
and that all their unit spheres satisfy $S(x) \in \mathcal{X}_{d-1}$.
The induction starts with $\mathcal{X}_{-1}=\{\}$. Having the class $\mathcal{X}_{d-1}$
of $(d-1)$-dimensional Dehn-Sommerville spaces topologically invariant immediately
bootstrap to see that also $d$-dimensional Dehn-Sommerville spaces have the property
that they are invariant under homeomorphisms. 

\satz{The edge refinement operation produces a homeomorphic graph and so
of its Whitney complex. If a graph is a $d$-manifold, then the edge refined
graph is a $d$-manifold. One can look at edge refinements as {\bf local
Barycentric refinements}. We also mentioned that Dehn-Sommerville spaces, a 
class of graphs generalizing spheres and like spheres forming a submonoid of
all complexes, are topological in nature. A homeomorphic sibling of a Dehn-Sommerville
space is Dehn-Sommerville.  }

\section{Fundamental group}

\paragraph{}
A {\bf closed curve} in a graph $G$ can be defined as a continuous map from a 
circular graph $C_n$ to $G$. This means that the vertices $x_0,x_1,\dots,x_n=x_0$ 
are mapped into vertices $y_0,y_1,\dots ,y_n=y_0$ such that either $y_i=y_{i+1}$ or
$(y_i,y_{i+1}) \in E$.  The {\bf fundamental group} of a graph equipped with a 
reference point $v$ is defined as the equivalence classes of 
closed curves in $G$ starting at $v$ modulo {\bf curve homotopy deformations}, 
where two curves are called {\bf curve homotopic}, if they can be morphed into each other by 
homotopy steps. 

\paragraph{}
A {\bf homotopy step} is an operation, where an edge of the path
attached to a triangle $t$ (an embedded complete graph $K_3$) is replaced with 
the two other sides or 
then reverses such a homotopy deformation and replaces two edges of the path in a triangle
with the other edge. If $f:G \to H$ is a continuous map, then a closed path maps either into 
a closed path or a point. The following result mirrors corresponding results in the continuum.
It is however a statement in finite topological spaces.

\paragraph{}
Replacing $1$-spheres by a $d$-sphere $S$ equipped with a base point, one can look at sphere 
embeddings continuous images of $S$ attached to a base point in $G$ and so look at 
homotopy groups $\pi_n(G)$. The sum of two such embeddings $S_1,S_2$ can be obtained by embedding
a wedge sum. We were once interested in graph complements of circular graphs \cite{GraphComplements}
because there, all higher dimensional wedge sums of spheres appear (at least homotopically equivalent)
as graphs complements of cyclic graphs. It would be nice if one could use this to compute
higher homotopy groups better but this has not worked yet. 

\satz{If $f: G \to H$ is a continuous map on graphs, it induces a 
group homomorphism $f_*: \pi_1(G) \to \pi_1(H)$ on the fundamental groups.}

\section{Euler characteristic}

\paragraph{}
The {\bf Euler characteristic} of a finite abstract simplicial complex $\mathcal{G}$ is defined as
$\chi(\mathcal{G}) = \sum_{x \in \mathcal{G}} \omega(x)$, where $\omega(x) = (-1)^{{\rm dim}(x)}$. 
The quantity can be seen in different ways. It is first of all a {\bf valuation},
meaning that it satisfies $\chi(\mathcal{G} \cup \mathcal{H}) = \chi(\mathcal{G}) + \chi(\mathcal{H}) -
\mathcal{G} \cap \mathcal{H}$. If $f_k(G)$ counts the number of $k$-dimensional simplices then also
$\chi(\mathcal{G}) = \sum_{k=0}^{\infty} (-1)^k f_k(G)$. If $G=(V,E)$ is a graph, its Euler
characteristic $\chi(G)$ is defined as the Euler characteristic of its Whitney complex $\mathcal{G}$. 
If $f:V \to \mathbb{R}$ is a locally injective function on vertices, 
then $i_f(v) = 1-\chi(S^-(v))$ is the Poincar\-Hopf index of $f$ 
at the vertex $v$. By induction, one can check the {\bf Euler-Poincar\'e formula} $\chi(G) = \sum_v i_f(v)$. 
Applying this to the graph $G_1$ of the Whitney simplicial complex $\mathcal{G}$ of a graph $G$
and using the function $f(x)= {\rm dim}(x)$ which is locally injective, one immediately can see that
the Euler characteristic of $G$ and $G_1$ are the same. The reason is that $i_f(v)= 1-\chi(S^-(v)) = 
1+(-1)^k = \omega(v)$ because in the case of the dimension functional, $S^-(x)$ is the boundary sphere 
complex of $x$ which has by the Euler-Gem formula the 
Euler characteristic $1+(-1)^k$ if $x$ has dimension $k$. The Euler characteristic therefore preserves
Barycentric refinements. One can see this also by explicitly writing down the linear map
transforming the {\bf $f$ vector} $(f_0,f_1,\dots,f_d)$ of $G$ to the $f$ vector of its 
Barycentric refinement. There is only one eigenvalue $1$ of this linear operator $T$ and 
the corresponding eigenvector of $T^*$ defines the Euler characteristic. 

\paragraph{}
One can see from the Poincar\'e-Hopf formula immediately that homotopy extensions and homotopy
reductions preserve the Euler characteristic: choosing a function $f$ which has the property that
it is maximal on the added vertex, we get $i_f(v) = 1-\chi(S^-(v))=0$ because
$S^-(v)$ is contractible and because recursively one sees that contractible graphs have Euler
characteristic $1$. One can also see that edge refinements in general preserve
the Euler characteristic: if $e=(a,b)$ is an edge, then the edge refinement replaces the 
join of $K_2$ with of $S(a) \cap S(b)$ with the join of the path graph $P_3$ with 
$S(a) \cap S(b)$. The operation just replaces a contractible part with an other contractible part
meaning that the Euler characteristic of that part does not change. Similarly, one can show
other operations like {\bf flip diagonal operations} on embedded kite graphs do not 
change the Euler characteristic. But flip diagonal operations does not preserve d-manifolds
in general. 

\paragraph{}
If simplices in $\mathcal{G}$ are equipped with an orientation 
\footnote{There does not need to be compatibility with intersecting simplices}
one can interpret an arbitrary function $f: \mathcal{G} \to R$
as a differential form. For $y \subset x$, define ${\rm sign}(y,x)=1$ if the orientation of 
$y$ matches the orientation of $x$ restricted to $y$, and $-1$ else. 
The {\bf exterior derivative} $df(x) = \sum_{y \subset x, |x|-|y|=1} {\rm sign}(y,x) f(y)$ 
is a $n \times n$ matrix if $\mathcal{G}$ has $n$ elements. Because if $z \subset y \subset x$
with $|z|=|y|-1$ satisfies $\sum_y {\rm sign}(z,y) {\rm sign}(y,x)=0$, the matrix $d$ satisfies
$d^2=0$ so that the {\bf Hodge Laplacian} $L= d d^* + d^*d$ is block diagonal with $f_k \times f_k$ block matrices 
$L_k$ leaving invariant the class of {\bf $k$-forms}, functions on $k$-dimensional simplices. 
The kernel of $L_k$ is called the {\bf $k$'th Betti number} of 
$\mathcal{G}$. By using the McKean-Singer symmetry that the non-zero eigenvalues of $L$ on even forms
agrees with the non-zero eigenvalues of $L$ on odd forms, one can see that the super trace
${\rm str}(A) = \sum_k (-1)^k A_{kk}$ has the property that $\chi(\mathcal{G}) = {\rm str}(e^{-t L})$
for any $t$. For $t=0$, one has ${\rm str}(1) = \sum_k (-1)^k f_k(G)$ and 
in the limit $t \to \infty$, where only the kernels of $L_k$ survives, one gets
$\chi(G) = \sum_k (-1)^k b_k$. The identity $\sum_k (-1)^k f_k(G) = \sum_k (-1)^k b_k(G)$ 
is called the {\bf Euler-Poincar\'e formula}.

\paragraph{}
The Betti numbers of the Barycentric refinement $\mathcal{G}_1$ are the same than the Betti 
numbers of $\mathcal{G}$. This can be seen as a consequence of the {\bf K\"unneth formula} 
which relates the Betti numbers of $H \cdot G$ with the Betti numbers of $H$ and the Betti numbers of 
$G$. The Stanley-Reisner product $H \cdot G$ is homotop to the Shannon product $H * G$ for which 
one can show the K\"unneth formula by taking the product of Harmonic functions $d^*f g$.
Homotopy deformations preserve the Betti numbers. 
If $\mathcal{H}= f(\mathcal{G})$ is a continuous image of $\mathcal{G}$ then 
$b_k(\mathcal{H})  \leq b_k(\mathcal{G})$. From these statements one can get immediately
that homeomorphic geometries have the same Betti numbers and so the same Euler 
characteristic. 

\satz{
Betti numbers, cohomology group, Euler characteristic are topological invariants.
Also the sphere spectrum $\bigcup_{x \in G} \chi(S(x))$ is a topological invariant.
The valuation property $\chi(U \cup V) = \chi(U) + \chi(V) - \chi(U \cap V)$
holds for Euler characteristic and all subsets $U,V$ of $G$. 
The Euler-Poincar\'e identity $\sum_k (-1)^k f_k(G) = \sum_k (-1)^k b_k(G)$ can 
be seen by heat deformation and using McKean-Singer symmetry. 
}

\section{Characteristics} 

\paragraph{}
Euler characteristic is the first of many {\bf higher characteristics}.
The next after Euler characteristic is {\bf Wu characteristic}. It is defined as
$\omega(\mathcal{G}) = \sum_{x \cap y \in \mathcal{G}} \omega(x) \omega(y)$.
Also all higher characteristics are invariant under Barycentric refinements.
For manifolds with boundary, it is $\chi(G) - \chi(\delta G)$ (see \cite{valuation}). 
It is a multi-linear valuation but not a valuation. It had been puzzling to 
us how the Wu characteristic behaves, even when looking at wedge sums. 
While the Euler characteristic is invariant under homotopy and so also under homeomorphisms
the Wu characteristic is only invariant under homeomorphisms. 

\paragraph{}
To see why topology is involved, we have to restate that the energy 
theorem tells $\chi(G) = \sum_{x,y \in G} g(x,y)$. There is a
quadratic identity to that $\omega(G) = \sum_{x,y \in G} \omega(x) \omega(y) g(x,y)^2$
(see \cite{GreenFunctionsEnergized}).  
Because $g(x,y) = \omega(x) \omega(y) \chi(U(x)\cap U(y))$ is 
expressed in terms of the topology. What happens is that
$f^{-1} (U(x) \cap U(y))$ is an open set with the same Euler characteristic. What happens is that
if $U,V$ are arbitrary open sets in the topology of $\mathcal{G}$, then 
$\omega(U \cup V) = \omega(U) + \omega(V)-\omega(U \cap V)$. 

\paragraph{}
The Euler characteristic $\chi(\mathcal{G})$ as a linear combination of basic
valuations $f_k(\mathcal{G})$ counting simplices. 
It satisfies the {\bf valuation formula} 
$\chi(\mathcal{G} \cup \mathcal{H}) = \chi(\mathcal{G}) + \chi(\mathcal{H}) - \chi(\mathcal{G} \cap \mathcal{H})$
if $\mathcal{G}$ and $\mathcal{H}$ are arbitrary simplicial complexes. This formula does not 
hold for the Wu characteristic. For example, if $\mathcal{G}$ is the {\bf Octahedron complex} and 
$\mathcal{H}$ is the {\bf circle complex} $C_4$, then 
$\omega(\mathcal{G}) = \chi(\mathcal{G})=2, \omega(\mathcal{H})=\chi(\mathcal{H})=0$. 
If we look at the {\bf wedge sum} $\mathcal{G} \wedge \mathcal{H}$ which is $\mathcal{G} \cup \mathcal{H}$
with a common $1$ point complex $K_1$. Now $\omega(K_1)=1$. We compute
$\omega(\mathcal{G} \cup \mathcal{H}) = 3$ so that obviously, the valuation formula does work as in the case of
Euler characteristic, where $\chi(\mathcal{G} \cup \mathcal{H})=1$. Now we know the solution to the 
puzzle: while $\mathcal{G},\mathcal{H}, \mathcal{G} \cap \mathcal{H}$ are all open sets by themselves,
in the topology of $\mathcal{G} \cup \mathcal{H}$, the complexes $\mathcal{G},\mathcal{H}$ are only closed in 
$\mathcal{G} \cup \mathcal{H}$ and no more open. We have however the identity
$\omega(U \cap V) = \omega(U) + \omega(V) - \omega(U \cap V)$ for open sets within the topological space
$\mathcal{O}$ in $\mathcal{X} = \mathcal{G} \cup \mathcal{H}$. 
If we look at the open balls $A=\mathcal{G} \setminus \overline{\{x\}}$,
$B=\mathcal{H} \setminus \overline{\{x\}}$ (they are open as an open set 
intersected with the complement of a closed set), and the open set $C=U(x)$ in $\mathcal{X}$. 
Now $A,B,C$ are open sets in $\mathcal{X}$ and 
we have $\omega(A \cup B \cup C) = \omega(A) 
+ \omega(B) +\omega(C) - \omega(A \cap B) - \omega(B \cap C) + \omega(A \cap B \cap C)$.

\paragraph{}
Consider the figure 8 graph $\mathcal{X}$ which is the wedge sum of two circular graphs 
$\mathcal{X} = \mathcal{G} \sup \mathcal{H} =C_4 \wedge C_4$. 
We can explain $\omega(\mathcal{X})=7$ by putting things together.
The valuation formula does not work for closed sets. For example, the following formula does not work: 
$7=\omega(\mathcal{X}) = \omega(\mathcal{G}) + \omega(\mathcal{H})
-\omega(\mathcal{G} \cap \mathcal{H}) = 0+0-1$. However, we can write $\mathcal{G}$ as a union
of three open sets $U,V,W$. Both $U,W$ are linear graph without boundary which have 
Wu characteristic $1$. The star graph without boundary has Wu characteristic $9$. 
The intersection between $U$ and $V$ has Wu characteristic $2$. So, we have
$\omega(\mathcal{X}) =\omega(U)+\omega(V)+\omega(W)-\omega(U \cap V)-\omega(V \cap W)
= 1+9+1-2-2=7$. 

\paragraph{}
We first used the old definition 
$\omega(U) = \sum_{x,y, x \cap y \neq \emptyset} \omega(x) \omega(y)$ for Wu characteristic
and not the correct definition $\omega(U) = \sum_{x,y, x \cap y \in U} \omega(x) \omega(y)$. 
There is no difference between the two definitions if we deal with simplicial complexes
which are closed sets. It matters however if we deal with open sets 
for example, take the two open sets $U = \{ (1,2,3),(1,2) \}$ and $V=\{ (1,2,3),(2,3) \}$,
which are stars in the complete complex $K_3$. Now $W=U \cap V = \{ (1,2,3) \}$ is the star of 
the facet $(1,2,3)$ which has $\omega(W)=1$. We have $X=U \cup V = \{ (1,2),(2,3),(1,2,3) \}$ 
with $\omega(x)=-1$. We have $\omega(U)=\omega(V)=0$ and $\omega(W)=1$, $\omega{X}=-1$. 
The identity $\omega(U) + \omega(V) - \omega(U \cap V) = \omega(U \cup V)$ is valid but
only because the simplices $x=(1,2)$ and $y=(2,3)$ were not allowed to ``interact". 
Their intersection was not in $U \cup V$. 

\paragraph{}
A convenient way to compute Wu characteristic therefore is to write the complex as a union
$\bigcup_j U_j$ of open sets, then use the {\bf inclusion exclusion property} 
$\omega(\mathcal{G}) = \sum_j \omega(U_j) - \sum_{i \cap j} \omega(U_i \cap U_j)
+ \sum_{i \cap j \cap k} \omega(U_i \cap U_j \cap U_k)$. This explains again the 
known fact that for $d$-manifolds $M$ we have $\omega(M) = \chi(M)$. 
Lets assume now that $G$ is a manifold graph with vertex set $V$. We can cover
$G$ with the open sets $U(v), v \in V$. Now use that $\omega(U(x)) = 1$ for any 
simplicial complex. So we have 
$\omega(G) = \sum_{x=(v) \in G} \omega(U(x))  - \sum_{x=(v,w) \in G} \omega(U(x)) + ... $
which is $\sum_x (-1)^{{\rm dim}(x)} \omega(x) = \chi(G)$. For manifolds $\omega(G)=\chi(G)$. 

\paragraph{}
In general we have the {\bf star formula}: 
$$ \omega(\mathcal{G}) = \sum_{x \in \mathcal{G}} \omega(x) \omega(U(x)) \;  $$
using the stars $U(x)$ of the simplex $x$.  We had previously proven the formula 
$\omega(\mathcal{G}) = \sum_{x,y \in \mathcal{G}} \omega(x) \omega(y) \chi(U(x) \cap U(y))^2$
in \cite{GreenFunctionsEnergized} we have here a sum over $\mathcal{G}$ and not a more
costly sum over pairs in $\mathcal{G}$.  We also see the {\bf ball formula} for the unit balls
$B(x)=\overline{U(x)}$ which is remarkable because there is no direct relation between
$\omega(B(x))$ and $\omega(U(x))$. The relation
$$ \omega(\mathcal{G}) = \sum_{x \in \mathcal{G}} \omega(x) \omega(B(x)) \;  $$
follows from $\sum_{x \in \mathcal{G}} \omega(x) \omega(S(x))= 0$. 
(See \cite{KnillEnergy2020} Corollary 6).

\satz{
Wu characteristic is no homotopy invariant but a topological invariants. 
The valuation property $\chi(U \cup V) = \chi(U) + \chi(V) - \chi(U \cap V)$
holds for all open sets,
where $\omega(U) = \sum_{x,y, x \cap y \in U} \omega(x) \omega(y)$ for a set of sets $U$
and not $\omega(U) = \sum_{x,y, x \cap y \neq \emptyset} \omega(x) \omega(y)$. 
There is the Gauss-Bonnet type formula 
$\omega(\mathcal{G}) = \sum_{x \in \mathcal{G}} \omega(x) \omega(U(x))$ complementing
$\chi(\mathcal{G}) = \sum_{x \in \mathcal{G}} \omega(x)$. This allows to compute the
Wu characteristic for larger spaces. There is also an energy 
theorem $\omega(G) = \sum_{x,y} g(x,y)$, where 
$g(x,y) = \omega(x) \omega(y) \omega(U(x) \cap U(y))$.  }

\section{Dynamics}

\paragraph{}
If $T$ is a continuous map from a finite topological space $O$ into itself
then every point is eventually periodic. Similarly, a simplicial map $T$ on
a simplicial complex.
The attractor of $T$ is a finite set and on every connected component of the 
attractor, one just cyclically permutes points.
One calls the forward attractor also the {\bf $\omega$-limit set}. 
In the case of a homeomorphism, there is also the {\bf $\alpha$ limit set} which 
is the $\omega$-limit set of the inverse map. 

\paragraph{}
The {\bf Lefschetz fixed point theorem} for graphs 
\cite{brouwergraph} tells that if $T$ is a graph endomorphism $T: G \to G$ 
then the sum of the indices of the fixed points agrees with the {\bf Lefschetz 
number} $\chi_T(mathcal{G})$, the super trace $\sum_{k} (-1)^k {\rm tr}(L|H^k)$ 
of $T$ induced each space of harmonic forms $H^k={\rm ker}(L_k)$.
This result formulated for graphs \cite{brouwergraph} 
obviously works for arbitrary simplicial complexes and a 
{\bf continuous map} $T: \mathcal{G} \to \mathcal{G}$. If $\mathcal{F}$ 
is the set of fixed points of $T$ and index $i_T(x) = \omega(x) {\rm sign}(T|x)$
with ${\rm sign}(T|x)$ being signature of the permutation of $T$ induced on $x$,
then the Lefschetz formula tells $\sum_{x \in \mathcal{F}} i_T(x) = \chi_T(\mathcal{G})$. 
The formula is easy to prove using the heat flow. The Koopman operator 
$U: f \to f(T)$ has as the super trace $\sum_{x \in \mathcal{F}} i_T(x)$. 
Applying the heat flow does not change the super trace of $e^{-t L} U$ 
as non-zero eigenvalues in the odd forms and even forms agree. In the limit
$t \to \infty$, only the map induced on the kernel survives and this is the 
Lefschetz number. A special case of the Lefschetz fixed point theorem is the 
Brouwer fixed point theorem which applies in the case when the complex has trivial 
cohomology. A special case is if $\mathcal{G}$ is contractible. An even 
more special case is if $\mathcal{G}$ is a $d$-ball. 

\paragraph{}
We can also start with an arbitrary finite topological space $\mathcal{O}$ 
and fix a pre-basis $\mathcal{B}$ which generates the topology. This defines a
{\bf nerve} simplicial complex $\mathcal{G}$ on $\mathcal{B}$, where the complex consists of all subset
of $\mathcal{B}$ which have a non-empty intersection. A continuous map $f$ on $\mathcal{O}$ 
now defines a continuous map on the simplicial complex $\mathcal{G}$ so that the 
{\bf Lefschetz fixed point theorem} applies. We can now define the cohomology of the topological space
(equipped with the base) as the cohomology of $\mathcal{G}$. The Lefschetz number of the super trace of 
the from $T$ induced map on the cohomology is then equal to the sum of the indices of fixed points of 
$f$ on $\mathcal{G}$. This means that there is an open set in $\mathcal{O}$ which is fixed.

\paragraph{}
A homeomorphism $T: G \to G$ of a finite geometry can be enhanced to a sequence
of homeomorphisms $T_n: G_n \to G_n$. As more iterations are needed, as more Barycentric 
refinements are required. For a homeomorphism $T$ this means
specifying a sequence of permutations $T_n: \mathcal{O}_n \to \mathcal{O}_n$ 
of the topologies of $\mathcal{G}_n$ and then require some compatibility. How well the 
map $T_n$ on $\mathcal{O}_n$ approximates the dynamics of $T_m$ on $\mathcal{O}_m$ determines
the {\bf amount of regularity of smoothness}.

\paragraph{}
The upgrade of a homeomorphism $T: \mathcal{G} \to \mathcal{G}$ to 
a {\bf stratified sequence of homeomorphisms} $T_n: \mathcal{G}_n \to \mathcal{G}_n$ is motivated by 
various similar constructions in mathematics, like computing with sequences of rational numbers with
a larger and larger number of digits in order to approximate real numbers, to do numerical computations 
of partial differential equations on sequences of meshes or the concept of inverse limit in 
constructions like p-adic integers or then martingales, where a stochastic process is observed on 
a sequence of adapted $\sigma$-algebras. 
Since every stochastic process given in the form of a sequence of IID 
random variables $X_n$ can be assigned a compact topological space $\Omega$ a continuous
function $f$ and a homeomorphism $T$ such that $X_n = f(T^n)$. The sigma-algebra
$\mathcal{A}_n$ generated by the random variables $X_1,\dots ,X_n$ is the Borel $\sigma$
algebra of a topological space $\mathcal{O}_n$ which is the product space $\Omega^n$. 
If $X$ has a finite set as range, then $\mathcal{O}_n$ is a finite topological space
and $\mathcal{A}_n$ is the Borel $\sigma$ algebra generated by $\mathcal{O}_n$. 

\satz{The Lefschetz fixed point theorem and so the Brouwer fixed point theorem 
naturally work for continuous maps on simplicial complexes. It even works for a finite
topological space when applied to the nerve of a pre-basis.
In order to study the dynamics of a homeomorphism, one has to pick a choice of 
concrete homeomorphisms on refinements. The length of the orbit which one 
wants to compute accurately determines how many Barycentric refinement lifts are needed.}

\section{Categorical}

\paragraph{}
Here are some general contemplations about the various categories: complexes, graphs and
topologies involved. Finite simplicial complexes form a category ${\rm Sim}$ with simplicial 
maps as morphisms. Finite graphs form a category ${\rm Gra}$ with graph homomorphisms as 
morphisms. Finite topological spaces ${\rm Top}$ form a category too where continuous maps are
the morphisms. We have the Whitney map from graphs to complexes, the Alexandroff map from 
complexes to topological spaces and the \v{C}ech map from topological spaces to the nerve
graph. These maps are not functors because the morphisms do not correspond directly. We can 
however enlarge the class of morphisms both on ${\rm Sim}$ as well as on ${\rm Gra}$ to have
morphisms. For example, in order to see the {\bf Whitney map} ${\rm Gra} \to {\rm Sim}$ as 
a {\bf functor} between categories one has to expand the possible morphisms on graphs 
allowing not only graph homomorphisms but maps from one graph to an other in which edges 
can collapse to vertices. In order to have a functor between ${\rm Sim}$ and ${\rm Top}$, 
we enlarge the class of simplicial maps and allow also continuous maps, still order preserving
but not mapping simplicial subcomplexes to simplicial subcomplexes necessarily (these are the
open maps). Also the map assigning to a simplicial complex a graph is a functor
again if one uses the larger class of morphisms. The composition of the functors 
from Graphs to Complexes back to Graphs is the Barycentric refinement. 
If ${\rm Sim}/~$ be the equivalence classes of complexes under 
Barycentric refinement and ${\rm Gra}/~$ the equivalence classes 
of graphs under Barycentric refinement. The Whitney map 
now identifies these two categories. A continuous map from some $G_n \to H$ could 
serve as the morphism. We chose however to make more assumptions and use the notion 
of having $\mathcal{H}$ a continuous image of $\mathcal{G}$ to define homeomorphisms. 

\paragraph{}
Some of the graph theory literature assumes graphs are
{\bf one-dimensional simplicial skeleton complexes}. While useful 
for some set-up's it is rather limiting as graphs are so much more than
one dimensional objects. The Whitney complex reflects rather
general topological spaces. There are other simplicial complexes associated to graphs 
of course, like the {\bf graphical matroid} {\bf other skeleton complexes} 
or the {\bf neighborhood complex}.

\paragraph{}
Topological graph theory looks at graphs embedded in two-dimensional manifolds. In that case, a
graph naturally naturally defines a cell complex in which the two-dimensional faces are the 
connected components of the complement of the embedded graph. This uses {\bf infinity} but it 
allows to deal with two-dimensional complexes which have the discrete topology of the underlying surface. 
On a surface of degree $g$ for example, the number of vertices $v$, the number of edges $e$ and the 
number of faces $f$ satisfies $v-e+f=2-2g$. Also in topological graph theory one can sometimes avoid infinity.
The notion of being planar for example is settled with Kuratowski's theorem completely within finite 
mathematics. That theorem uses {\bf homeomorphism} in the narrow sense as homeomorphic as one-dimensional 
simplicial complexes. 

\paragraph{}
Discrete CW complexes extend simplicial complexes. 
A {\bf discrete CW structure} can be introduced within combinatorics,
once one has defined what a $k$-sphere is:
Start building up the geometry $G_0=\{\}$ and successively attach $k$-balls 
(called {\bf cells} or {\bf handles}) to already existing $(k-1)$-spheres. 
We can for example attach a $0$-ball (called a vertex) 
to a $-1$ sphere (the empty graph). Once the 0-dimensional part
is built, we can attach $1$-balls (called edges) to $0$-spheres (2 disjoint points). 
Then one can start adding 2-dimensional balls (faces) to $1$-spheres. 
For example, one can add triangles $2$-simplices to the
a triangular circle. Obviously, every finite abstract simplicial complex 
is also an abstract CW complex.  
While simplicial complexes are natural and given just as they are, a CW complex comes 
with a ``timeline" of how the structure has been built up. 
We can so build also {\bf multi-graphs} or {\bf quivers} or more
general complexes called {\bf $\delta$-sets}, which are simplicial complexes in 
which simplices can occur with multiplicities. 
Adding a bit more structure produces a subclass of $\delta$-sets 
called {\bf simplicial sets}. 
In category theory this is known as a {\bf pre-sheaf} on the simplex category. 
If on a $\delta$-set boundary maps are defined, one has a cohomology 
like on simplicial complexes. 

\paragraph{}
We have more recently also looked at {\bf quivers}, graphs where self-loops and 
multiple connections can happen. In that case, one can naturally attach {\bf $\delta$-sets}
to a quiver. $\delta$-sets generalize simplicial complexes in that different sets can appear
multiple times and where boundary maps are specified.
The category of $\delta$ sets generalize the category of {\bf simplicial sets}.
The later are $\delta$ sets with more structure attached. $\delta$ sets (and so also simplicial sets) have
a cohomology attached. One can now ask, what natural topologies can be associated to a 
quiver. We have not yet investigated this. One possibility would be to do this on the sets
$x$ in the $\delta$-complex and take the basis $U(x)$. But now, the space is not even $T_0$ 
any more as different points can have the same minimal open sets. One can no more distinguish
points by their minimal open neighborhoods. 

\satz{
Simplicial complexes, finite simple graphs and topological spaces can not be directly 
linked with their traditional morphisms. But topology glues them together if we look at continuous
maps as morphism. Whether we talk about a graph with 
continuous maps on them, simplicial complexes with continuous maps on them or finite topological 
spaces with continuous maps on them, we always can switch to the other two pictures.  
It requires however to change already what we mean by morphisms. While traditionally, these
three categories use different notation and jargon, topology unifies them nicely and allow us
to work on finite geometries using a trinity of interpretations. Non-standard analysis links
this radically finite geometry with rather arbitrary compact topological spaces. 
}

\section{Ringed spaces}

\paragraph{}
As in {\bf commutative algebra} approaches to geometry, one can use the notion of 
{\bf ringed space}. Attach a ring $F(U)$ to ever open set and call it a {\bf section} of $U$.
Given restriction maps produces a {\bf pre-sheaf}.
The usual way to rephrase this is that this is a contra-variant functor from the 
category of open sets with inclusion morphisms to the category of rings. 
To get a sheaf, we need existence (gluing) and uniqueness (locality) properties:
Gluing is related to existence 
$h(x)| U(x) \cap U(y) = h(y)|U(x) \cap U(y)$ then there {\bf exists} a $h$ with $h|U(x)=h(x)$.
Locality relates to uniqueness because $h \in R(x) = k \in R(x)$ for all $x$, then 
we have the same $h=k$.

\paragraph{}
A simple case is to the sheaf of ring-valued continuous 
functions on open sets such that if $V \subset U$, the restriction of 
$F(U)$ to $F(V)$ is a ring homomorphism.  In our case, where the topology of a simplicial complex 
with the topology, we deal with a {\bf ringed simplicial complex}. 
The {\bf section} $F(x)=F(U(x))$ is in this context the {\bf stalk} of $x$ and its
elements are the {\bf germs}. Given a commutative local ring $R$ and any 
function $h: \mathcal{G} \to R$ defines already a {\bf locally ringed sheaf}. But things
can be much more general. The restriction maps from $U(x) \to U(y)$ if $x \subset y$
do not have to be the obvious ones. Actually, any ring valued matrix $r(x,y)$ 
can serve as a transition map $F(x) \to F(y)$ if $x \subset y$. 
The pre-sheave condition now means for $x \subset y \subset z$ that the cocycle condition
$r(y,z) r(x,y) = r(x,z)$ holds and especially that $r(y,x) r(x,y) = r(x,x)$.

\paragraph{}
Classically, when looking at general topological spaces, 
the {\bf stalk} $F(U(x))$ at some point $x$ is the {\bf direct 
limit} $F(U)$ over all the open sets $U$ containing $x$. 
In the finite topology case, the stalk at $x$ is $F(U(x))$, 
which is just a ring attached to the star $U(x)$. 
A {\bf locally ringed space} is a ringed space in which every stalk is a 
local ring, meaning that it has a unique maximal ideal.
An example of a ringed space are differential forms. 
If an orientation is fixed on each simplex $x \in \mathcal{G}$,
then these are just the functions from $\mathcal{G}$ to the ring $R$.

\paragraph{}
Finite geometries also allow to use the frame work of {\bf schemes} in an elementary frame work. 
A ringed space by definition attaches to every open set a ring and its spectrum, 
the set of {\bf prime ideals} in the ring. If $U$ is an open finite set, then the space $\mathcal{O}_X(U)$ of
functions on $U$ have as prime ideals the functions which vanish at some simplex $x$. 
The spectrum therefore is just the set of simplices in $U$. 
The theory as developed for general ringed spaces can be taken over word for word. The frame work can be useful
also in combinatorics. For example, we can take $\mathbb{C}$-valued functions and require that 
the restriction maps from  ring $F(U)$ to the ring $F(V)$ if $V \subset U$
is not the obvious ones. Locally ringed topological space can have
global properties are not necessarily the obvious ones: 
going around a closed loop for example can produce a non-trivial map.

\paragraph{}
The lack of linear structures prevents having constructs like tangent spaces in the discrete.
However, we have attached to each simplex a unit sphere $S(x)$ and so a 
{\bf sphere bundle}. What we can do in general is to have transition maps on 
spheres $S(x)$ coming from positive dimensional simplices. These transitions
tell what happens if one looks at $S(x)$ as part of $S(v)$ or $S(w)$ if 
$v,w \subset x$. Going from a pre-sheaf to a sheaf means to have transition
maps on positive dimensional simplices. 

\satz{
Having a topology on a complex allows to use sheaf theoretical
concepts on finite spaces. In finite topological spaces, 
the ring $R(x)$ attached to a star is called a stalk
and its elements are the germs. 
}

\section{Morse extensions}

\paragraph{}
A function $f: \mathcal{G} \to R$ to a totally ordered space $R$ like $\mathbb{R}$ of $\mathbb{Z}$ 
is called a {\bf Morse function} if it is {\bf locally injective},
meaning that $f(x) \neq f(y)$ if $x \subset y$ or $y \subset x$ and $S^-(x) = \{ y \in S(x), f(y) < f(x) \}$
is a $(k-1)$-sphere for some $k \geq 0$ or then contractible. In the former case, we have
added a {\bf critical point}, in the later case a {\bf regular point}.
If $x$ is a critical point, the integer $k \geq 0$ is called the {\bf Morse index}
of the point $x$. If $x_k$ is a fixed enumeration of points such that $f(x_k) \geq f(x_l)$ if
$k \geq l$, then $\mathcal{G}_n = \{ y, f(y)<f(x_n) \}$ is a Morse build-up
an $\chi(\mathcal{G}_{n+1})=\chi(\mathcal{G}_{n}) + \chi(B^+(x_n)) -\chi(S^-(x_n))$
$= \chi(\mathcal{G}_n) + i_f(x_n)$ so that the {\bf Poincar\'e-Hopf formula} 
$\chi(G) = \sum_x i_f(x)$  \cite{poincarehopf} holds. 
This formula holds for any locally injective function $f$ but for Morse functions
the Poincar\'e-Hopf index $i_f(x) \in \{-1,1\}$. We see that the existence of a Morse function
implies that $\mathcal{G}$ can be seen as a {\bf CW complex} in which successively $k$-balls 
\footnote{also called handles} are attached to $(k-1)$-spheres in the previous step. 

\paragraph{}
The Morse build-up of a graph are not homeomorphism even if we look at the step when a 
regular point is added. One can see this already from the fact that the dimension can increase
without adding a critical point. A continuous map can not increase dimension. 
The contraction process is however a continuous process.
The stable unit sphere $S_f^-(x)$ is then a subgraph of $S(x)$ and so a closed set. 
We can state this all in other words and say that a real valued function on a simplicial 
complex needs not to be continuous but that every function on the vertex set of a graph
is continuous in the topology of the Barycentric refinement. 

\paragraph{}
It would not be useful to enforce continuity because $\{ f(x)<c \}$ can be a single point
which is neither open nor closed in the topology of $\mathcal{G}$. 
When we look however at the situation on the graph level with a function 
$f: V(G) \to R$, then the sets $G_n= \{v, f(v) \leq v_n \}$ 
generate graphs which are closed sets in $\mathcal{G}1$. In this sense
any function $f: V \to R$ on a graph $G$ is continuous in the topology of $\mathcal{G}$ 
while a function on a-priori given simplicial complex $\mathcal{G}$ is hardly ever
continuous in that topology. A function $\mathcal{G} \to R$ becomes
only continuous if we look at it as a function of the graph $G_1$ and so using the topology 
of $\mathcal{G}_1$. 

\satz{
Morse theory is a concrete way to build up an abstract finite 
CW-complex using a Morse function as a guidance. 
While the sets $\mathcal{G}_n$ in a Morse build-up of a complex 
are neither open or closed, the topology of their graphs make them 
topological spaces in a refinement. Functions on simplicial complexes become naturally 
continuous when considered in the topology of the Barycentric refinement. }

\section{Calculus}

\paragraph{}
A function $f: \mathcal{G} \to R$ can be interpreted as a differential form. Similarly, 
if $G$ is a graph, we look at functions $f: \mathcal{G} \to R$. 
When restricted to $k$-dimensional simplices, one gets {\bf $k$-form}.
Calculus can be studied on arbitrary Barycentric refinement levels. Provided that orientations
are fixed on $\mathcal{G}$ any scalar function is just a differential form.
In order to define level surfaces, we need only the very mild assumption that 
functions are locally injective, meaning that adjacent vertices take different values. 
Lest look in this section at graphs. In case of a simplicial complex $\mathcal{G}$, 
look at the graph $G_1$ in which the vertex set is $\mathcal{G}$ and where two are connected
if one is included in the other. 

\paragraph{}
A locally injective function $f: G \to R$ can be lift to a function on $G_1$ 
by assigning to a simplex $x$ the average of $f$ over the vertices in $x$.  
We can also take a function  $f: V(G) \to R$ and distribute its values 
$f(v)$ equally to all points $x$ in $U(v)$. This produces a new 
function $G_1 \to R$. What might happen under such a refinement of a locally injective
function, that it is no more locally injective. We can fix this by looking at the lexicographic
order of the pair $(f,{\rm dim})$. Because in a Barycentric refinement, 
${\rm dim}(x) \neq {\rm dim}(y)$ if $x,y$ are connected in $G_1$, this is a well defined
ordering. We can now use $f$ on $G_1$ to build again a level surface.

\paragraph{}
Given a function on the vertex set of $G_1$ we can move the content from vertices
in $G_1$ which are sets of vertices in $G$ to vertices of
$G$ by equally distributing the value $f(x)$ to to all vertices $v \in x$. 
With $f(x) = \omega(x)$, this produces $f(v) = \kappa(v)$, 
where $\kappa(x)$ is the {\bf curvature} \cite{Levitt1992}. The consequence
$\sum_x \omega(x) = \sum_v \kappa(v)$ is the {\bf Gauss-Bonnet theorem.} 
See \cite{cherngaussbonnet,knillcalculus}. 

\paragraph{}
If $f$ is a function on a discrete $d$-manifold $G$ which is locally injective. 
Then the {\bf level surface} $U=\{ f = c \}$ generated by the set of simplices 
on which $f$ changes sign. See \cite{KnillSard}. 
In the topology of the complex, this is an open set. 
However, its graph defines a discrete $(d-1)$ manifold if it is not empty. This manifold
now carries a topology again. Note that as a set $M$ in $\mathcal{G}$, 
the set $U$ is always open because if $x$ is in $M$ then every set $y$ containing $x$
is in $M$. We can still make $f$ locally injective by replacing
$f$ with $f(x)+\epsilon {\rm dim}(x)$ for $\epsilon$ small. This allows us to extend
$f$ to $G_1$ in a determined way. We can now look at a {\bf variety} 
$\{f_1=0,\dots,f_m=0\}$ for $m$ locally injective functions 
$f_1, \dots, f_m$ in the $n$'th Barycentric refinement 
$G_n$ as the set of simplices, where all the lifted functions 
of $f_j$ change sign. 

\paragraph{}
We have seen that using the ``dimension trick" providing a lexicographic
order of the function on a higher level, we can lift any function on a graph uniquely to 
Barycentric refinement where we can again define level sets. 
This level set is a graph where two simplices are connected if one is contained 
in the other. Remarkably, by the {\bf discrete Sard theorem},
we never run into singularities.
We expect if the functions $f_k$ are lifted nicely to a Barycentric refinement, 
then the corresponding manifold is homeomorphic to $S$. We expect that there
could be surprises if we take a situation from the continuum, where 
$\{f_1=0,\dots,f_m=0\}$ is a classical variety which is not a manifold. In that
case, there could be surprises near singularities. The topology depends on 
the choice of the Barycentric refinement. 

\satz{If $G$ is a $n$-manifold and $m \leq n$ locally injective functions
are given on $G$, then the ``variety" $\{f_1=0,\dots,f_m=0\}$ is a well defined
graph again. It is either empty or a $(n-k)$-manifold $S$. }

\section{Interaction energy}

\paragraph{}
Given a simplicial complex $\mathcal{G}$ with $n$ elements $x$ and any
$n \times n$ matrix taking values in some ring $R$, we can define the 
{\bf internal energy} of a subset $A \subset \mathcal{G}$ as 
$$ \omega(A) = \sum_{x,y, x \cap y \in A} h(x,y) \; . $$
The matrix $h$ does not have to be symmetric. We can think of $h(x,y)$
also as a current from $x$ to $y$ and $\omega(A)$ as the total current or traffic 
flowing overall through $A$. Now look at the matrix 
$$  g(x,y) = \sum_{x,y} \omega(x) \omega(y) \omega(U(x) \cap U(y)) \;.  $$
This matrix gives a {\bf potential energy} between the simplices $x$ and $y$. 
Unlike $h$, the matrix $g$ is always symmetric. 
The {\bf energy theorem} tells $\omega(\mathcal{G}) = \sum_{x,y} g(x,y)$. 
This result assures that the total potential energy of $\mathcal{G}$ is the total internal
energy. This theorem generalizes an energy theorem proven before, see 
\cite{KnillEnergy2020,GreenFunctionsEnergized,EnergizedSimplicialComplexes2}. 
For example, if $h(x,y)$ is diagonal with $h(x,x) = \omega(x)$, then $\omega(A) = \chi(A)$
is the Euler characteristic. 
In that case the matrix $g$ is the inverse of the operator
$L(x,y) = \chi(\overline{x} \cap \overline{y})$, where $\overline{x}$ is the closure of $\{x\}$,
a simplicial complex. The frame work also captures {\bf energized simplicial complexes}
where $h(x,x) = h(x)$ and $h(x,y) = 0$ for $x \neq y$. 

\paragraph{}
To prove the more general energy theorem, 
note that the map $h \to g$ is linear and that both the energy and the total
sum are both linear expressions. We only need to verify the statement therefore
for the matrix $h$ satisfying $h(x_0,y_0)=1$ and $h(x,y)=0$ else for some fixed simplices $x_0,y_0$. 
These $n^2$ basis elements are fixed-points of the linear map $T(h)=g$ and
the energy relation $\sum_{x,y, x=x \cap y \in G} h(x,y) = \sum_{x,y} g(x,y)$ holds.
By linearity, the relation holds then for all $h$.
To verify the statement for a basis element, 
note that the left hand side is $1$ if $x_0$ and $y_0$ intersect and $0$ else. 
The right hand side is $\omega(x) \omega(y) \omega(U(x) \cap U(y))$ which is non-zero only if $x_0$
and $y_0$ intersect and both $x_0,y_0$ are contained in $U(x)$ and $U(y)$.
This means that the union $x_0 \cup y_0$ is contained in $U(x) \cap U(y)$.
This means that both $x$ and $y$ have to contain $x_0 \cup y_0$ . 
This means that the simplex $x \cap y$ has to contain the simplex $x_0 \cup y_0$ and so $x_0 \cap y_0$. 
But $\sum_{x_0 \subset x,y_0 \subset y, x_0 \cap y_0 \subset x \cap y} \omega(x) \omega(y)
=\sum_{x_0 \subset x} \omega(x) \sum{y_0 \subset y} \omega(y) =1*1 = 1$
because the Euler characteristic of a simplex is $1$. 

\paragraph{}
We also have as before a curvature relation 
$\kappa(x) = \sum_{y} g(x,y) = \omega(x) g(x,x) = \omega(x) \chi(U(x))$ and so 
$\sum_{x} \kappa(x) = \chi(G)$ which is a {\bf Gauss-Bonnet relation}. It can also 
be thought of as a {\bf Poincar\'e-Hopf theorem} for the locally injective function $f(x) = -{\rm dim}(x)$
because then, the atom $U(x)$ is the stable sphere $S^-(x) = \{ y, f(y) < f(x) \}$. 
Still, since the {\bf dimension function} $f$ is not so well visible, it is good to think of $\kappa$
not as an {\bf index} but as a {\bf curvature}. 
The internal energy of the ``atom" $U(x)$ is up to a sign a curvature. Summing up 
the curvature gives the total energy. We can also think of the relation $\sum_x \omega(x) g(x,x)$
as the {\bf super trace} ${\rm str}(g)$ of the matrix $g$. In the context of simplicial complexes, 
the notion of super trace is natural since ${\rm str}(1) = \chi(G)$ is the Euler characteristic
and because of the {\bf McKean-Singer formula} ${\rm str}(e^{-t L}) = \chi(G)$. 
\cite{knillmckeansinger,McKeanSinger}. 

\paragraph{}
Actually, we would like to announce here already that arbitrary {\bf tensor energy theorems} hold. Let 
$h(x_1,\dots,x_m)$ be arbitrary ring-valued functions of $m$ variables. We already had for
$$ \omega_2(A) = \sum_{x,y,x \cap y  \in A} h(x,y)  \;  $$
that  $g_2(x,y) = \omega(x) \omega(y) \omega_2(U(x) \cap U(y))$ satisfies $\sum_{x,y} g_2(x,y) = \omega_2(G)$. 
If  $\mathcal{G}$ is a finite abstract simplicial complex and $h(x,y,z)$ arbitrary $R$ valued function.
For a subset $A \subset \mathcal{G}$, define the {\bf internal cubic energy}
$$ \omega_3(A) = \sum_{x,y,z,x \cap y \cap z \in A} h(x,y,z)  \; . $$
Now define $g_3(x,y,z) = \omega(x) \omega(y) \omega(z)  \omega_3(U(x) \cap U(y) \cap U(z))$
and think about it as the {\bf potential energy} of the triple. 
Then the total potential energy agrees with the total internal energy 
$$   \sum_{x,y,z} g_3(x,y,z) = \omega_3(G) \; .  $$
This works also with more interaction like quartic
$$ \omega_4(A) = \sum_{x,y,z,w, x \cap y \cap z \cap w \in A} h(x,y,z,w)  \; . $$
For $g_4(x,y,z,w) = \omega(x) \omega(y) \omega(z) \omega(w) \omega(U(x) \cap U(y) \cap U(z) \cap U(w))$,
the total potential energy agrees with the total internal energy 
$$   \sum_{x,y,z,z} g_4(x,y,z,w) = \omega_4(G) \; .  $$

\satz{
Given an interaction transfer rule between $m$ intersecting simplices
in a simplicial complex $\mathcal{G}$, 
we can assign internal energies $k$ tuples of sets. The total
energy can be expressed also as the sum of all potential energies.
The internal $m$-energy of a set $\omega(A) = \sum_{\bigcap_j x_j \in A} h(x_1,\dots,x_m)$ 
is associated to closed sets, the potential 
energy uses open sets $g_m(x_1,\dots,x_k) = \prod_j \omega(x_j) \omega(\bigcap_{j=1}^k U(x_j)) $. 
The energy theorem assures that the total internal energy is the total
potential energy $\omega_m(\mathcal{G}) = \sum_{x_1,\dots,x_k} g(x_1,\dots, x_k)$.
In the case $m=1$ and $h(x)=\omega(x)$, where the total energy is the Euler
characteristic, this is a Gauss-Bonnet theorem 
$\chi(G)=\sum_x \omega(x) \omega(U(x))$, where $U(x)$ is the smallest open set
containing $x$. In the case $h(x_1,\dots,x_m) = \prod_{k=1}^m \omega(x_k)$, the energy 
$\omega(\mathcal{G})$ is the $m$-th characteristic, a topological invariant for the complex. 
Unlike for $m=1$, which gave the Euler characteristic, we have no homotopy invariant 
however for $m>1$. 
}

\section{Remarks}

\paragraph{}
Inspired by Alexandroff and Zariski, we have revisited here at the finite
topology on a simplicial complex $\mathcal{G}$ defined by stars and especially for
simplicial complexes coming from a finite simple graph $G=(V,E)$. 
\footnote{Since open sets are still quite local, the 
drawbacks of Zariski topology appear not really relevant.}
The analogy to algebra is that the {\bf vertices} of the graph play the role 
of the {\bf maximal ideals} and that the {\bf simplices} play the role of
{\bf prime ideals}. A basis for the topology is the set of stars of a simplex, 
the set of simplices which contain $x$, 
as well as an added empty set. Unlike any topology on the vertex set $V$ like the
one given by distance which would render the graph completely disconnected,
our topology honors connectivity and dimension. It shares the non-Hausdorff property
with the Zariski topology on prime ideals of a commutative ring.
The closed sets in our graph topology are exactly the simplicial complexes of subgraphs.
They play the role of algebraic subsets of a variety in algebraic geometry. 
The closure of a set $A$ of simplices is the smallest abstract simplicial
complex which contains the set $A$.

\paragraph{}
One of the motivations to look at finite topologies is that there
are finite simple spaces $G,H$ for which the geometric realization $|G|,|H|$ 
are homeomorphic but which are not homeomorphic in the finite topology. 
Open sets $U(x),U(y)$ be entangled in a
complicated way in the finite topology if the dimension is large. But these entanglements
are not always visible when looking at the topology induced from Euclidean distance in
geometric realizations. The topology of the Euclidean realization is not sophisticated
enough. In other words, there are triangulations of topological manifolds which have
the manifold as a geometric realizations but which are not discrete manifolds as 
defined here. The definition of homeomorphism is motivated by the notion of {\bf piecewise linear}
map in topology. A map $f: |\mathcal{G}| \to |\mathcal{H}|$ between geometric realizations
of simplicial complexes is PL, if there is a piecewise linear map between Barycentric
refinements $|\mathcal{G}_n| \to |\mathcal{H}_m|$. 
A {\bf PL} homeomorphism is then a simplicial map such that there is a 
homeomorphism $|\mathcal{G}_n| \to |\mathcal{H}_m|$ for
some refinement. This looks equivalent to what we do here in the finite but leaves finite 
mathematics. As we are not interested in infinity here, we do not bother showing the equivalence.

\paragraph{}
An other motivation has been to answer the question why {\bf higher characteristics} like {\bf Wu
characteristic} $\omega(\mathcal{G}) = \sum_{x,y, x \cap y \in \mathcal{G}}  \omega(x) \omega(y)$
is a {\bf topological notion} while {\bf Euler characteristic} 
$\chi(\mathcal{G}) = \sum_x \omega(x)$ is more. We will write about this more in the future but
one of the key facts is that $\omega(U \cup V) = \omega(U) + \omega(V) - \omega(U \cap V)$ 
holds for open sets but not for closed sets or sets which are neither closed nor open in general. 
For Euler characteristic this {\bf valuation formula} holds for all subsets $U,V$ of $\mathcal{G}$. 
Euler characteristic does super count simplices, while Wu characteristic does super count intersecting
simplices. This interacting points should not be separable. In indeed, if an intersecting pair $x,y$  
is in the intersection $U \cap V$ of two open sets it must be in both sets $U$ and $V$. An other
mystery which still needs more investigation is the notion of {\bf analytic torsion}
$A(\mathcal{G}) = \prod_k {\rm Det}(L_k)^{k (-1)^{k+1}}$,
where $L_k$ are the blocks of the Hodge Laplacian $L=D^2=(d+d^*)^2$ of the Whitney
complex and ${\rm Det}$ is the pseudo determinant. We wrote this as a super determinant
${\rm SDet}(D) = \prod_k {\rm Det}(D_k)^{(-1)^k}$ of the Dirac operator $D$ with Dirac blocks
$D_k=d_k^* d_k$. Analytic torsion can be defined for any simplicial 
complex but so far it has been accessible only in 2 cases: the first is when $\mathcal{G}$ is homotopic
to $1$. In that case $A(\mathcal{G}) =|V|$ where $V=\bigcup_{x \in G} x$. 
The second case was when $\mathcal{G}=|V| |V'|$ for even dimensional spheres and 
$\mathcal{G} = |V|/|V'|$ for odd dimensional spheres, where $V'$ is the number of maximal simplices
in $\mathcal{G}$. But analytic torsion is {\bf not} a topological invariant. Even for manifolds like
a torus, the analytic torsion changes if we make a Barycentric refinement. As already in the etymology 
of the name, a non-trivial fundamental group makes the functional $A(\mathcal{G})$ more complicated as
we can get torsion terms along non-contractible closed loops. Still it is not only that. If we make
homotopy extensions of a sphere which are not homeomorphisms, the torsion formula for the sphere in 
terms of $|V|$ and $|V'|$ disappears. 

\paragraph{}
Finite and so Alexandroff topologies can be interesting from a physics point of view. 
First of all, there is local interaction of 
simplices which are contained in each other. We can not separate such points using
open sets. The fact that we have smallest non-empty open {\bf Planck units} $U(x)$ 
is some sort a {\bf space quantization} or {\bf atoms of space}. If we look at a manifold with a very fine 
triangulation, both the lack of the Hausdorff topology and the Alexandroff feature 
are hardly visible. The situation is also present in floating point arithmetic, 
when a computer deals with small numbers. Every point has a smallest neighborhood 
which can no more be resolved. We also can not separate two 
points which are too close from each other even so they are different.
Points which are identified by the equal-tolerance parameter define
the {\bf machine graph}. With a machine precision $\log_{10}(2^{52}) \sim 15.65$
the distance below which two numbers are identified is about 
$2^{-46} = 2^{-52-1+7} = 1.42109 \cdot 10^{-14}$.

\paragraph{}
The definition of homeomorphism is motivated by the fact that if we have two finite
topological spaces coming from a finite abstract simplicial complex 
and a continuous surjective map $f: X \to Y$ and a continuous surjective map $g: Y \to X$
and the unit spheres $S(x)$ have homeomorphic pre-images and unit balls of locally maximal
simplices have balls as images, then $X,Y$ are homeomorphic.  
\footnote{We wonder whether it is true in general: 
does already the existence of two continuous surjective maps $f:X \to Y, g:Y \to X$ 
force $X,Y$ to be homeomorphic. Using the axiom of choice, one can invert the 
surjections and have injections, showing with Cantor-Schroeder-Bernstein that the cardinalities 
are the same so that there is a bijection between the topologies. 
In the non-Alexandroff case, where points can be written as 
intersections of open sets, this should then give a homeomorphisms.  }

\paragraph{}
The given definition of {\bf homeomorphism within finite topology} 
has shifted a bit while writing down this text. We first tried to avoid 
unit spheres $S(x)$ and balls but failed to prove some results. 
Using unit sphere $S(y)$ in the definition allows the use of {\bf induction}. For the unit balls
of maximal simplices, there is no interesting topology as they are balls and we just require that
the inverse image of such a unit ball is a ball (which is precisely defined). It is postulated that
all $d$-dimensional balls are homeomorphic but it could also be proven from the definition. 
It also seems to be necessary to ask such a requirement. We need it for
example that if $\mathcal{G}$ is a manifold and $\mathcal{H}$ is homeomorphic, then $\mathcal{H}$
is a manifold. We can imagine continuous surjective maps going both ways which collapse substantial
parts of the topology somewhere in the interior of the manifold so that the inverse image of $S(y)$
could become a complicated object in $\mathcal{G}$.

\paragraph{}
Looking at stars $U(x)=W^+(x)$ and cores $V(x) = \overline{ \{x\}}=W^-(x)$ 
of simplices is also motivated by {\bf connection calculus}. 
We have seen for example that the {\bf Green function matrix}
$g(x,y)=\omega(x) \omega(y) \chi(U(x) \cap U(y))$
is always the inverse to the matrix $L(x,y) = \chi(V(x) \cap V(y))$, where
$\chi(A) = \sum_{x \in A} \omega(y)$ is the Euler characteristic of an arbitrary
subset of $\mathcal{G}$ and $\omega(x) = (-1)^{{\rm dim}(x)}$. While $V(x) \cap V(y)$
which is always a simplicial complex, are closed, the sets $U=U(x) \cap U(y)$ are open. 
The Euler characteristic of the closure $B(x)=\overline{U}$ is in general different from the Euler
characteristic of $U$. Actually $\chi(B(x)) = \chi(U(x)) + \chi(S(x))$ where $S(x)$ is the 
boundary of $U$.  \footnote{This phenomenon prevented us for some time to find
the Green star formula for the Green function.}

\paragraph{}
Both the star $U(x)$ and the core $V(x)=\overline{\{x\}} = W^-(x)$ can be seen as 
{\bf measurable sets} in the graph. 
If we close a topology $\mathcal{O}$ under complements and countable intersections and
unions, we get the {\bf Borel $\sigma$ algebra} $\mathcal{A}$ as usual.
This set $\mathcal{A}$ still does not cover all subsets of $\mathcal{G}$
if $G$ has dimension $2$ or more. The reason is that if $\{e=(a,b)\}$ is a single
non-maximal simplex, then it is neither open nor closed. The smallest open set
containing it is $U(e)$, the smallest closed set containing it is the
simplicial complex $\{ e,\{a\},\{b\} \}$.
We can look for probability theory on the graph and look for example at the
measure $f_k(A)/f_k(G)$ counting the fraction of $k$-dimensional simplices in $A$.
For probability theory on finite set, see \cite{Nelson}.

\paragraph{}
In the literature, one often a ``graph" as a topological space that is obtained as a
geometric realization as a {\bf one-dimensional simplicial complex}. 
A more topological approach is to look at the geometric realization of its Whitney 
complex in which all the complete subgraphs $K_{n+1}$ are realized as simplices. One 
can then also look at other simplicial complexes attached to a graph, similarly
as one can attach other topologies to $\mathbb{R}^n$. The notion of homeomorphism could be
extended to such cases too. There are simplicial complexes and so graphs that are not 
homeomorphic but which have homeomorphic realizations. An example is a double suspension of a
rational homology sphere. It is not topologically equivalent to a sphere in our sense but in 
a geometric realization, it is by the {\bf double suspension theorem}.
While also in the discrete, any manifold that is a suspension of a manifold 
must be a sphere, in the discrete, a suspension of a non-manifold is by definition
not a discrete manifold.

\paragraph{}
We have searched for notions of homeomorphism within finite combinatorics for a while
like \cite{KnillTopology}, where we looked at \v{C}ech type notions like the {\bf nerve}
of an open cover and asked that two homeomoprhic graphs have isomorphic nerves. 
In 2016 we experimented (motivated by the Zarisiki topology) 
with the concept of having the closed subgraphs play the role of closed sets. 
We have not taken it too seriously: do we want to work with topological spaces that are non-Hausdorff? 
We decided now it is better to work with the topology generated by the star basis. 
When reviewing the Lusternik-Schnirelmann category, where open covers play a role, the concept
fits better with what one does in the continuum. Since 2016, we have also realized
more the importance of {\bf stars} $U(x)=U(x)$ and {\bf cores} $\overline{x} = W^-(x) =\overline{\{x\}}$ 
of simplices as they form a hyperbolic structure and because 
$g(x,y) = \omega(x) \omega(y) \chi(U(x) \cap U(y))$ and
$L(x,y) = \chi(W^-(x) \cap W^-(y))$. In the current notation, we would write this {\bf Green-Star formula} as
$g(x,y) = \omega(x) \omega(y) \chi(U(x) \cap U(y))$. The fact that $U(x) \cap U(y)$ can be topologically 
quite complicated even so both sets $U(x),U(y)$ are smallest open sets and have contractible 
closures, the intersection $U(x) \cap U(y)$ can be rather complicated. The local smallest Planck units $U(x)$
in a complex can be entangeled in a complicated way. Also, the simplicial complex belonging to the 
closure of $U(x) \cap U(y)$ can be topologically very different from $U(x) \cap U(y)$!
\footnote{This was a reason to drove us almost insane in 2016 while looking for the Green star formula as the
formula with the closure $\overline{U(x) \cap U(y)}$ worked in most cases and especially small complexes, 
but that it had rare failures.}

\paragraph{}
There are many interesting open questions and many opportunities for experimentation or
further explorations. We can ask for example about the fraction 
$|\mathcal{O}(G)|/2^{|\mathcal{G}(G)|}$ telling us in a graph what fraction of subsets 
of the simplicial complex are open sets. As the number of open sets and closed sets agree, 
this is equivalent to count the number of subgraphs of a given graph. 
Numerical computations become quickly too hard to do. We have $\phi(C_4)=48/256$ and 
$\phi(C_5)=124/1024$ and $\phi(C_6)=323/4096$ and $\phi(K_1)=1$, $\phi(K_2)=5/8$,
$\phi(K_3)=19/128$ and $\phi(K_4)=167/32768$. 

\paragraph{}
Let us add a comment on the literature. The paper \cite{Alexandroff1937}, 
was dedicated to Emmy Noether, assumes that the space is locally finite spaces which under
a global compactness assumption means finiteness. Alexandroff already identifies discrete
$T_0$ spaces with partially ordered sets and identifies closes sets as simplicial complexes.
He notes that if in an Alexandroff space, two smallest open sets $U(x)=U(y)$ agree, then $x=y$. 
The reason is that then $x \in U(y)$ and so $y \subset x$ and $y \in U(x)$ and so $x \subset y$.
He calls simplicial complexes {\bf vollst\"andige Mengensystems} = complete set systems. The completion 
of an arbitrary set of sets is the closure. Interestingly, the void $\{ \}$ = empty set is 
not considered of this type, even so today we consider this to be a simplicial complex. It is technically
a finite set of sets which is closed under the operation of taking finite non-empty subsets. Also
in topology, the empty set is a closed set as in any topological space we consider $\emptyset,
X$ to be clopen = closed and open and define connectedness as the property that the only clopen 
sets are $\emptyset$ and $X$. The modern point of view is to see the void $0$ as the $(-1)$ dimensional sphere.
Alexandroff also notes that pre-base of stars centered at $0$-dimensional simplices define what we 
would call today a \v{C}ech graph. In our terminology we would say that 
every locally finite simplicial complex coming from a graph $G$ has a \v{C}ech cover (the pre-base) whose 
graph is $G$. It has also a \v{C}ech cover (coming from the base) which is the Barycentric refinement.
Alexandroff also notes that a continuous map can be lifted to Barycentric 
refinements. He however refers to the geometric realization as a polyhedron in order to define
something analog to homeomorphism. Alexandroff also reformulates the construction of the cohomology ring 
following Alexander-\v{C}ech and Whitney. 

\paragraph{}
Finite topologies spaces were picked up again as such in the 1960ies like \cite{Stong1965}.
Strong for example shows that for any finite topological space, there is a unique minimal base.
Connectedness and path connectedness are equivalent.
That a continuous self-map $f$ on a finite topological space that is either injective or surjective
must be a homeomorphism. Strong shows already that on finite topological spaces, continuity 
is equivalent with simplicial map: $x \subset y$ if and only if $f(x) \subset f(y)$.  
Strong  equips the space $H^G$ of continuous maps $G \to H$ with the compact-open topology. 
It is ordered with $g \leq f$ if for all points $g(x) \leq f(x)$. If $f \leq g$ then $f,g$
are homotopic. The connectivity components of $H^G$ are the homotopy classes of maps from $G$ to $H$. 
\cite{McCord1966} starts with {\it Finite topological spaces have more interesting topological
properties than one might suspect at first.} Indeed, McCord points out that for every finite
topological space, there is a finite simplicial complex which is a weak homotopy equivalent
in the sense that the induced maps on all homotopy groups are isomorphisms (meaning for $\pi_0$
which is not equipped with a group structure, that the number of connected components are the same). 
McCord is known also for a version of the nerve theorem stating that the homotopy type of a nice topological 
space is encoded in the \v{C}ech nerve of a nice open cover. This certainly applies for finite topological
spaces and the cover coming from minimal open sets. What is needed for example is that the intersection of two
such sets is either contractible or empty. 
The \v{C}ech nerve of a cover has been introduced by Alexandroff. 

\pagebreak

\section{Code} 

\paragraph{}
The following few Mathematica lines allow to compute the topology of a complex or the 
topology of the Whitney complex of a graph. We see that the number of topologies
on the cyclic graph $C_n$ is the Lucas number $L(2n)$. We then display the 
code for Figure 1. 

\begin{tiny}
\lstset{language=Mathematica} \lstset{frameround=fttt}
\begin{lstlisting}[frame=single]
Closure[A_]:=If[A=={},{},Delete[Union[Sort[Flatten[Map[Subsets,A],1]]],1]];
Whitney[s_]:=If[Length[EdgeList[s]]==0,Map[{#}&,VertexList[s]],
   Map[Sort,Sort[Closure[FindClique[s,Infinity,All]]]]];
UU[G_,x_]:=Module[{U={}},Do[If[SubsetQ[G[[k]],x],
   U=Append[U,G[[k]]]],{k,Length[G]}];U];
Basis[G_]:=Table[UU[G,G[[k]]],{k,Length[G]}];
SubBasis[G_]:=Module[{V=Union[Flatten[G]]},
   Table[UU[G,{V[[k]]}],{k,Length[V]}]];
UnitSpheres[G_]:=Module[{B=Basis[G]},
      Table[Complement[Closure[B[[k]]],B[[k]]],{k,Length[B]}]];
UnitBalls[G_]:=Map[Closure,Basis[G]];
Cl[U_,A_]:=Module[{V=U},Do[V=Union[Append[V,
   Union[V[[k]],A[[l]]]]],{k,Length[V]},{l,Length[A]}];V];
Topology[G_]:=Module[{V=B=Basis[G]},
   Do[V=Cl[V,B],{Length[Union[Flatten[G]]]}];Append[V,{}]];
GraphBasis[s_]:=Basis[Whitney[s]];
GraphTopology[s_]:=Topology[Whitney[s]];
Nullity[Q_]:=Length[NullSpace[Q]];
Fvector[G_]:=Delete[BinCounts[Map[Length,G]],1];
ToGraph[G_] :=Module[{n=Length[G],v,e,s},v=Range[n];e={};
   Do[If[(SubsetQ[G[[k]],G[[l]]]|| SubsetQ[G[[l]],G[[k]]]) &&   
   Not[G[[k]]==G[[l]]],e=Append[e,v[[k]]->v[[l]]]],
   {k,n},{l,k+1,n}]; s=UndirectedGraph[Graph[v,e]]];
BarycentricGraph[s_]:=ToGraph[Whitney[s]]; 
BarycentricComplex[G_]:=Whitney[ToGraph[s]]; 
w[x_]:=-(-1)^Length[x]; 
Wu1[A_]:=Total[Map[w,A]]; Chi=Wu1; 
Wu2[A_]:=Module[{a=Length[A]},Sum[x=A[[k]]; Sum[y=A[[l]];
   If[MemberQ[A,Intersection[x,y]],1,0]*w[x]*w[y],{l,a}],{k,a}]]; Wu=Wu2; 
Wu3[A_]:=Module[{a=Length[A]},Sum[x=A[[k]]; Sum[y=A[[l]]; Sum[z=A[[o]];
   If[MemberQ[A,Intersection[x,y,z]],1,0]*w[x]*w[y]*w[z],{o,a}],{l,a}],{k,a}]];
FastChi[A_]:=Module[{UU=Basis[A]},Sum[w[A[[k]]]*Chi[UU[[k]]],{k,Length[A]}]];
FastWu[A_]:=Module[{UU=Basis[A]}, Sum[w[A[[k]]]*Wu[UU[[k]]],{k,Length[A]}]];
FastWu3[A_]:=Module[{UU=Basis[A]},Sum[w[A[[k]]]*Wu3[UU[[k]]],{k,Length[A]}]];
Suspension[G_]:=Module[{q=Max[Flatten[G]]+1,n=Length[G]},
  Closure[Union[Table[Append[G[[k]],q],{k,n}],Table[Append[G[[k]],q+1],{k,n}]]]];
JoinAddition[A_,B_]:=Module[{q=Max[Flatten[A]],Q,G=A},Q=Table[B[[k]]+q,{k,Length[B]}];
  Do[G=Append[G,Union[A[[a]],Q[[b]]]],{a,Length[A]},{b,Length[Q]}];G=Union[G,Q];
  If[A=={},G=B]; If[B=={},G=A]; G];
DoubleSuspension[G_]:=Suspension[Suspension[G]];
WuBetti[G_]:=Module[{Cohomology2,n,n2,G2,ll,ln,dd1,dd2,LL2,L2,dd,br,D2,DD},
  n=Length[G]; length[x_]:=Length[x[[1]]]+Length[x[[2]]]; G2={}; IS=Intersection;
  Do[If[Length[IS[G[[k]],G[[l]]]]>0,G2=Append[G2,{G[[k]],G[[l]]}]],{k,n},{l,n}];
  n2=Length[G2]; G2=Sort[G2,length[#1]<length[#2] & ]; ll = Map[length,G2];
  ln=Union[ll];br=Prepend[Table[Max[Flatten[Position[ll,ln[[k]]]]],{k,Length[ln]}],0];
  derivative1[{x_,y_}]:=Table[{Sort[Delete[x,k]],y},{k,Length[x]}];
  dd1=Table[0,{n2},{n2}]; Do[u=derivative1[G2[[m]]]; If[Length[u]>0,
  Do[r=Position[G2,u[[k]]]; If[r!={},dd1[[m,r[[1,1]]]]=(-1)^k],{k,Length[u]}]],{m,n2}];
  derivative2[{x_,y_}] :=Table[{x,Sort[Delete[y,k]]},{k,Length[y]}];
  dd2 = Table[0,{n2},{n2}]; Do[u = derivative2[G2[[m]]]; If[Length[u]>0,
  Do[r=Position[G2,u[[k]]]; If[r!={}, dd2[[m,r[[1,1]]]]=(-1)^(Length[G2[[m,1]]] + k)],
   {k, Length[u]}]], {m,n2}]; dd = dd1 + dd2; D2=dd+Transpose[dd]; L2 =D2.D2;
  LL2=Table[Table[L2[[br[[k]]+i,br[[k]] + j]],{i,br[[k+1]]-br[[k]]},
   {j, br[[k + 1]]-br[[k]]}], {k,Length[br]-1}];
  Cohomology2=Map[NullSpace,LL2]; Map[Length,Cohomology2]];
\end{lstlisting}
\end{tiny} 

\paragraph{}
As example computations, we compute the number of elements in the
topology of a circle $C_n$ where $|\mathcal{G}|=2n$. One can show by induction in $n$
that the number of open sets in the topology is $L(2n)$, where $L(n)$
is the {\bf Lucas number} defined by $L(0)=2,L(1)=1,L(2)=3$ and 
the recursion $L(n+1)=L(n)+L(n-1)$ is the {\bf Fibonacci recursion}. (The only 
difference is that for the Lucas numbers, the entry $L(0)=2$, while for the 
Fibonacci numbers, the entry $F(0)=1$ is assumed.)  

\begin{tiny}
\lstset{language=Mathematica} \lstset{frameround=fttt}
\begin{lstlisting}[frame=single]
Table[Length[GraphTopology[CycleGraph[k]]],{k,4,7}]
Table[LucasL[2n],{n,4,7}]
\end{lstlisting}
\end{tiny}

\paragraph{}
Here we compute the Euler characteristic and Wu characteristic of 
star graphs:

\begin{tiny}
\lstset{language=Mathematica} \lstset{frameround=fttt}
\begin{lstlisting}[frame=single]
Table[s=StarGraph[k];{Chi[Whitney[s]],Wu[Whitney[s]]},{k,3,10}]
\end{lstlisting}
\end{tiny}

\paragraph{}
Here we compute the Wu characteristic of the basis of a random graph

\begin{tiny}
\lstset{language=Mathematica} \lstset{frameround=fttt}
\begin{lstlisting}[frame=single]
s=RandomGraph[{15,54}]; Map[Wu,GraphBasis[s]]
\end{lstlisting}
\end{tiny}

\paragraph{}
This is the code for Figure 1

\begin{tiny}
\lstset{language=Mathematica} \lstset{frameround=fttt}
\begin{lstlisting}[frame=single]
e={1->2,2->3,3->1,3->4,3->6,3->8,8->9,8->10}; V=ViewVertical;
s=UndirectedGraph[Graph[e]];BG=BarycentricGraph;
A=GraphPlot3D[s,ViewPoint->{1,-3,-1},V->{1,-1,-0.3}];
B=GraphPlot3D[BG[BG[s]],ViewPoint->{0,-2.5,-2},V->{1,0,0}];
S=GraphicsRow[{A,B}]; Export["figure1.pdf",S,"PDF"]; Show[S]
\end{lstlisting}
\end{tiny}

\paragraph{}
We illustrate Gauss-Bonnet 
$$  \omega(\mathcal{G}) = \sum_{x \in \mathcal{G}} \omega(x) \omega(U(x)) 
    = \sum_{x \in \mathcal{G}} w(x) \omega(B(x)) $$
and $\sum_{x  \in \mathcal{G}} \omega(x) \omega(S(x)) = 0$. 
These formulas hold for any simplicial complex $\mathcal{G}$. We have seen an 
analog formula $\chi(\mathcal{G}) =  \sum_{x \in \mathcal{G}} \omega(x) \chi(U(x))$ 
for Euler characteristic $\omega_1(\mathcal{G}) = \chi(\mathcal{G})$ before. But 
it holds for any higher characteristic $\omega_m(\mathcal{G})$.  

\begin{tiny}
\lstset{language=Mathematica} \lstset{frameround=fttt}
\begin{lstlisting}[frame=single]
s=RandomGraph[{44, 220}]; G=Whitney[s]; 
{ Timing[FastWu[G]], Timing[Wu[G]]} 
U=Basis[G]; S=UnitSpheres[G]; B=UnitBalls[G];
{Wu[G],Sum[w[G[[k]]]*Wu[U[[k]]],{k,Length[G]}],
       Sum[w[G[[k]]]*Wu[S[[k]]],{k,Length[G]}],
       Sum[w[G[[k]]]*Wu[B[[k]]],{k,Length[G]}]}
\end{lstlisting}
\end{tiny}

\paragraph{}
For small complexes, the direct Wu computation is faster. 
But already if $G$ has several hundred entries, the fast
Wu computation is faster. In the first example, where the 
complex had 59 elements, the direct computation was 4 times faster.
In the third of the following cases the fast Wu computation took 23 seconds 
while the Wu computation took 90 sections. The complex had 1355 elements. 

\begin{tiny}
\lstset{language=Mathematica} \lstset{frameround=fttt}
\begin{lstlisting}[frame=single]
s=RandomGraph[{14, 30}]; G=Whitney[s];
{ Timing[FastWu[G]], Timing[Wu[G]]}

s=RandomGraph[{44, 220}]; G=Whitney[s];
{ Timing[FastWu[G]], Timing[Wu[G]]} 

s=RandomGraph[{54, 420}]; G=Whitney[s];
{ Timing[FastWu[G]], Timing[Wu[G]]} 
\end{lstlisting}
\end{tiny}

\paragraph{}
In general, we measure  $\omega(U(x)) + \chi(B(x)) -\chi(S(x)) \geq 0$. 
This is equivalent to $\omega(U(x)) \geq \chi(U(x))$. This is something,
we have not been able to explain yet: 

\begin{tiny}
\lstset{language=Mathematica} \lstset{frameround=fttt}
\begin{lstlisting}[frame=single]
s=RandomGraph[{14, 30}]; G=Whitney[s];
U=Basis[G]; S=UnitSpheres[G]; B=UnitBalls[G];
Map[Chi, U] - Map[Chi, B] + Map[Chi, S]
Map[Wu, U] - Map[Chi, U]
\end{lstlisting}
\end{tiny}

\paragraph{}
We also have an {\bf energy theorem for Wu characteristic}
$$  \omega(G) = \sum_{x,y} g_2(x,y) \; ,   $$
where $g_2(x,y) = \omega(x) \omega(y) \omega(U(x) \cap U(y))$
is a Green function matrix. Unlike in the case of Euler characteristic,
where $g_1(x,y) = \omega(x) \omega(y) \chi(U(x) \cap U(y))$ was
unimodular, the matrix $g_2$ is no more unimodular in general. The 
determinant is in general not $1$. For manifolds, we compute it here for 
a 3-sphere, a double suspension of a cyclic graph $C_4$. 

\begin{tiny}
\lstset{language=Mathematica} \lstset{frameround=fttt}
\begin{lstlisting}[frame=single]
s=RandomGraph[{20, 40}]; G=Whitney[s]; n=Length[G]; U=Basis[G]; 
g2=Table[w[G[[k]]]*w[G[[l]]]*Wu[Intersection[U[[k]],U[[l]]]],{k,n},{l,n}];
Print[{Wu[G],Total[Flatten[g2]],Sum[w[G[[k]]]*g2[[k,k]],{k,n}]}]; 
Print[Det[g2]]
\end{lstlisting}
\end{tiny}

\paragraph{}
There would be a lot more to explore. We can look for example at the 
{\bf sphere Green matrix}
$$  s_2(x,y) = \omega(x) \omega(y) \omega(S(x) \cap S(y)) $$
and compare it with the {\bf ball Green matrix} 
$$  b_2(x,y) = \omega(x) \omega(y) \omega(B(x) \cap B(y)) $$
and the {\bf star Green matrix}
$$  g_2(x,y) = \omega(x) \omega(y) \omega(U(x) \cap U(y)) $$
for which we see that the trace and the total sum of all entries and the 
determinant are all zero. Is there some significance to the nullities we
see in the ball or sphere Green function entries. If Wu is replaced with Euler
we get Green function matrices $s_1,b_,g_1$. 

\begin{tiny}
\lstset{language=Mathematica} \lstset{frameround=fttt}
\begin{lstlisting}[frame=single]
s=RandomGraph[{20,40}]; 
G=Whitney[s];n=Length[G];U=Basis[G];S=UnitSpheres[G];B=UnitBalls[G]; 
s1=Table[w[G[[k]]]*w[G[[l]]]*Chi[Intersection[S[[k]],S[[l]]]],{k,n},{l,n}];
s2=Table[w[G[[k]]]*w[G[[l]]]* Wu[Intersection[S[[k]],S[[l]]]],{k,n},{l,n}];
Print[{Chi[G],Total[Flatten[s1]],Sum[w[G[[k]]]*s1[[k,k]],{k,n}],Det[s1]}];
Print[{Wu[G], Total[Flatten[s2]],Sum[w[G[[k]]]*s2[[k,k]],{k,n}],Det[s2]}];

b1=Table[w[G[[k]]]*w[G[[l]]]*Chi[Intersection[B[[k]],B[[l]]]],{k,n},{l,n}];
b2=Table[w[G[[k]]]*w[G[[l]]]* Wu[Intersection[B[[k]],B[[l]]]],{k,n},{l,n}];
Print[{Chi[G],Total[Flatten[b1]],Sum[w[G[[k]]]*b1[[k,k]],{k,n}],Det[b1]}];
Print[{Wu[G], Total[Flatten[b2]],Sum[w[G[[k]]]*b2[[k,k]],{k,n}],Det[b2]}];

g1=Table[w[G[[k]]]*w[G[[l]]]*Chi[Intersection[U[[k]],U[[l]]]],{k,n},{l,n}];
g2=Table[w[G[[k]]]*w[G[[l]]]* Wu[Intersection[U[[k]],U[[l]]]],{k,n},{l,n}];
Print[{Chi[G],Total[Flatten[g1]],Sum[w[G[[k]]]*g1[[k,k]],{k,n}],Det[g1]}];
Print[{Wu[G], Total[Flatten[g2]],Sum[w[G[[k]]]*g2[[k,k]],{k,n}],Det[g2]}];a

Print[{Nullity[g1],Nullity[b1],Nullity[s1],Nullity[g2],Nullity[b2],Nullity[s2]}] 
\end{lstlisting}
\end{tiny}

\paragraph{}
Gauss-Bonnet formulas for higher characteristic allows to compute higher invariants
more quickly. It still needs time. The {\bf homology 3 sphere} is implemented with a simplicial
complex of 392 simplices. \cite{BjoernerLutz}. It has Euler and Wu characteristic $0$. Its suspension has Euler 
and Wu characteristic $2$.  The double suspension again has Euler and Wu characteristic 0.
We know by the {\bf double suspension theorem} that the double suspension of the homology sphere
has a geometric realization that is homeomorphic to a 5 sphere. We unfortunately can not compute
the Wu cohomology yet as the complex is too large. We suspect that Wu cohomology can distinguish
the double suspension of the homology sphere from the 5 sphere. The following computation still needs
a few minutes to compute, even so we use the Gauss-Bonnet version for Wu characteristic. Already 
in the case of the Suspension of the homology sphere, the fast Wu computation is faster. For small
complexes, the direct computation is faster because the fast Wu procedure requires to pre-compute a
basis of the topology. 

\begin{tiny}
\lstset{language=Mathematica} \lstset{frameround=fttt}
\begin{lstlisting}[frame=single]
onesphere=Whitney[CycleGraph[4]];
moebius = Whitney[GraphComplement[CycleGraph[7]]]; 
cylinder=Whitney[UndirectedGraph[Graph[
   {1->2,2->3,3->4,4->1,5->6,6->7,7->8,8->5,1->5,5->2,2->6,6->3,3->7,7->4,4->8,8->1}]]];
twosphere=Suspension[onesphere]; 
threesphere=DoubleSuspension[onesphere];
homologyB={{1,2,4,9},{1,2,4,15},{1,2,6,14},{1,2,6,15},{1,2,9,14},{1,3,4,12},
{1,3,4,15},{1,3,7,10},{1,3,7,12},{1,3,10,15},{1,4,9,12},{1,5,6,13},{1,5,6,14},
{1,5,8,11},{1,5,8,13},{1,5,11,14},{1,6,13,15},{1,7,8,10},{1,7,8,11},{1,7,11,12},
{1,8,10,13},{1,9,11,12},{1,9,11,14},{1,10,13,15},{2,3,5,10},{2,3,5,11},{2,3,7,10},
{2,3,7,13},{2,3,11,13},{2,4,9,13},{2,4,11,13},{2,4,11,15},{2,5,8,11},{2,5,8,12},
{2,5,10,12},{2,6,10,12},{2,6,10,14},{2,6,12,15},{2,7,9,13},{2,7,9,14},{2,7,10,14},
{2,8,11,15},{2,8,12,15},{3,4,5,14},{3,4,5,15},{3,4,12,14},{3,5,10,15},{3,5,11,14},
{3,7,12,13},{3,11,13,14},{3,12,13,14},{4,5,6,7},{4,5,6,14},{4,5,7,15},{4,6,7,11},
{4,6,10,11},{4,6,10,14},{4,7,11,15},{4,8,9,12},{4,8,9,13},{4,8,10,13},{4,8,10,14},
{4,8,12,14},{4,10,11,13},{5,6,7,13},{5,7,9,13},{5,7,9,15},{5,8,9,12},{5,8,9,13},
{5,9,10,12},{5,9,10,15},{6,7,11,12},{6,7,12,13},{6,10,11,12},{6,12,13,15},{7,8,10,14},
{7,8,11,15},{7,8,14,15},{7,9,14,15},{8,12,14,15},{9,10,11,12},{9,10,11,16},{9,10,15,16},
{9,11,14,16},{9,14,15,16},{10,11,13,16},{10,13,15,16},{11,13,14,16},{12,13,14,15},
{13,14,15,16}}; homologysphere=Closure[homologyB];
Print[FastWu[homologysphere]];
Print[FastWu[Suspension[homologysphere]]]; 
Print[FastWu[DoubleSuspension[homologysphere]]];
\end{lstlisting}
\end{tiny}

\paragraph{}
Here are some computations of Wu Betti numbers. We could not yet complete
the computation of the Wu cohomology for the double suspension of the 
homnology sphere. 

\begin{tiny}
\lstset{language=Mathematica} \lstset{frameround=fttt}
\begin{lstlisting}[frame=single]
WuBetti[onesphere]       (* {0,1,1}               *)
WuBetti[twosphere]       (* {0,0,1,0,1}           *)
WuBetti[threesphere]     (* {0,0,0,1,0,0,1}       *)
WuBetti[moebius]         (* {0,0,0,0,0}           *)
WuBetti[cylinder]        (* {0,0,1,1,0}           *)
WuBetti[homologysphere]  (* not yet able to compute *)
WuBetti[DoubleSuspension[homologysphere]] (* dito   *)
\end{lstlisting}
\end{tiny}

\paragraph{}
{\bf Example 1}: For the smallest positive dimensional example $G=K_2$, 
we have the simplicial complex 
$\mathcal{G}=\{ (1),(2),(1,2) \}$. The basis has three elements and
consists of  \\
$\{ \{(1),(1,2)\}, \{ (2),(1,2) \},\{ (1,2) \} \}$. 
The corresponding unit spheres are  \\
$\{ \{ (2) \}, \{ (1) \}, \{ (1),(2) \} \}$. 
The topology has $5$ elements \\
$\{ \{ \}, \{(1),(1,2)\}, \{ (2),(1,2) \},\{ (1,2) \}, \mathcal{G} \}$.
Any path graph is homeomorphic to this graph.  \\
A homotopy reduction $f: K_2 \to K_1$ given by $1 -> 1, 2 -> 1$ 
is continuous. The topology in $K_1$ is $\{ \{ \}, (1) \}$ and the 
inverse of every of the open sets is an open set in $K_2$. There
is no continuous surjective map from $K_1$ to $K_2$. Any Barycentric
refinement of $K_1$ is $K_1$ and a map on finite spaces can not 
increase cardinality. Also, any continuous map can only decrease 
or preserve dimension.

\paragraph{}
{\bf Example 2}: All cyclic graphs $C_n$ with $n \geq 4$ are homeomorphic
but $C_n$ is not homeomorphic to a path graph $P_m$. There is no homeomorphism
as there would have to be a continuous surjective map
$f: C_n \to P_m$. This is not possible because there are two unit spheres $S(x)$ in $P_m$ for which
the inverse image has $2$ elements (a 0-sphere). As a 0-sphere $S_0$  is not homeomorphic
to a 1-point graph (there is not even a surjective map from $K_1$ to $S_0$).
One can also see that $C_n$ is not homeomorphic to $P_m$ because the Euler characteristic does not
match. One can also see it from the fact that the fundamental groups do not match. 

\paragraph{}
{\bf Example 3}: the definition of homeomorphism in finite spaces can be used to produce
constructive verifications that two spaces are homeomorphic or not. 
To do so, cover both spaces with balls which intersect in balls
then try to match the balls up. Obviously, if the maximal dimension of the two spaces is 
different they can not be homeomorphic. Let us assume that we have two complexes $\mathcal{G}$
and $\mathcal{H}$ which are homeomorphic and both have maximal dimension $d$,
then the number of connected components of $d$-dimensional maximal balls must be the same in 
both. Two star graphs $S(n)$ and $S(m)$ with different number of rays can not 
be homeomorphic for example. The star graph $S(n)$ has $n$ different open
one-dimensional balls, while $S(m)$ has $m$ different connectivity components. 

\section{Questions}

\paragraph{}
The definition of ``homeomorphism" proposed here seems have all the properties we want:
it has invariants like Euler characteristic, Wu characteristic, {\bf Lusternik-Schnirelmann category}
(the minimal number of contractible sets which cover the space), 
Betti numbers, Wu Betti numbers,  cup length, Lebesgue dimension, connectivity type, separation
properties or being a manifold are the same for homeomorphic graphs. 

\paragraph{}
There are other notions which are not topological invariants like the number of 
$k$-dimensional simplices $f_k$, the eigenvalues of some Hodge Laplacian, curvature,
inductive dimension, average simplex cardinality, Dehn-Sommerville invariants for non-manifolds, 
the Fermi characteristic $\phi(\mathcal{G}) = \prod_x \omega(x)$ which agrees with the 
determinant of the connection Laplacian ${\rm det}(L)$ \cite{KnillEnergy2020}. 
For Dehn-Sommerville, especially related to Gauss-Bonnet curvatures
\cite{cherngaussbonnet,valuation,dehnsommervillegaussbonnet}. We have shown for example that
the Dehn-Sommerville property is invariant under edge refinement, the join operation and
Barycentric refinements. Also Poincar\'e-Hopf 
\cite{poincarehopf,parametrizedpoincarehopf,PoincareHopfVectorFields,MorePoincareHopf}
can be reformulated more conveniently in a topological frame work. 

\paragraph{}
{\bf A)} One thing we could not explore yet whether the {\bf relative Wu characteristic} 
\cite{valuation,CohomologyWu} $\omega(\mathcal{G},\mathcal{H})$ for a subcomplex
$\mathcal{H}$ defined as $\sum_{x \in \mathcal{G}, y \in \mathcal{H},x \cap y \neq \emptyset} 
\omega(x) \omega(y)$ is depends on the topology of $\mathcal{H}$ and on the 
embedding in $\mathcal{G}$. The Wu characteristic 
$\omega(\mathcal{G}) = \sum_{x \cap y \neq \emptyset} \omega(x) \omega(y)$ 
itself is a topological invariant.  
We would also like to know to compute the Wu characteristic in a classical manner
without triangulation. Obviously, just looking at the structure of open covers 
does not work. What matters also are the local dimensions, the dimensions of the 
covers. For example, if we glue two manifolds along a k-dimensionial
part, then the dimension of this connection matters. If two graphs
$C_4$ are glued along a point (one calls this a {\bf wedge sum}), 
we get the figure 8 graph $G$ with $\omega(G)=7$. 
If we glue it along an edge, we get a {\bf digital figure 8 graph} $H$
and $\omega(H)=5$. These two graphs are not homeomorphic because there
are unit sphere which are not homeomorphic. The vertex degrees of 
$G$ are $2$ or $4$ while the vertex degree of $H$ are $2$ or $3$. 

\paragraph{}
In the context of calculus, there are questions about the minimal number of 
{\bf critical points} of locally injective functions. In general, a critical point is a point $x$
for which $S^{-}_f(x)  = \{ y \in S(x), f(y) < f(x) \}$ is non-contractible. Points for
which the Poincar\'e-Hopf index $i_f(x) = 1-\chi(S^{-}_f(x))\neq 0$ are critical points but
there can also be critical points of index $0$. The Lusternik-Schnirelman inequality 
assures ${\rm cup}(\mathcal{G})+1 \leq {\rm cat}(\mathcal{G}) \leq {\rm cri}(\mathcal{G})$ 
where ${\rm cup}$ is the {\bf cup length} (a homotopy and so topological invariant) and 
where ${\rm cri}$ counts the minimal number of critical points, which a locally injective
function can have. The Morse inequalities count the minimal number of critical points of 
a Morse function can have. Also here, one can ask whether this number is a topological 
invariant. More generally one can ask whether the numbers $c_k$ counting the minimal number
of Morse critical points of index $k$ is topological. 
The {\bf Morse inequalities} produce the 
general bound $\sum_k (-1)^k b_k \leq \sum_k (-1)^k v_k$. 

\paragraph{}
{\bf B)} As the main focus of this note was a definition of homeomorphism, 
it would be good to know more about redundancies in the definition. 
We have played with various versions of the definition. 

\begin{itemize}
\item We tried first not to make any requirement
about the unit ball and only have the condition for unit spheres.
\item An other modification would be to avoid talking about unit spheres and balls and 
ask that every unit ball $B(x)$ is homeomorphic to $f^{-1} B(x)$. (We have asked this only 
for locally maximal simplices where the unit ball is a ball)
This implies that the boundary $S(x)$ is homeomorphic
to the boundary of the inverse $f^{-1} B(x)$. It would also imply the for a maximal simplex
where $B(x)$ is a ball, the inverse $f^{-1} B(x)$ is a ball. We did not want to use this as
 a definition however because $B(x)$ is the same so that we have have no induction to work with. It still 
make the definition local. One could then postulate that all balls of the same dimension
are homeomorphic. But that is less elegant.
\item We also tried to play with the requirement that $f^{-1}(B(x))$ is contractible (which is weaker
than requiring it to be a ball) but we had difficulty from this to establish for example that 
if two spaces are homeomorphic and one is a manifold, the other must be a manifold.
\end{itemize} 

{\bf C)} The property that $\mathcal{H}=f(\mathcal{G}_n)$ is a continuous image was defined
as the property that any unit sphere $S(x)$ is homeomorphic to $f^{-1}(S(x))$ for all $x \in \mathcal{H}$
and that the inverse of the unit ball $B(x)$ is a ball in $\mathcal{G}_n$. 
Is this already is enough to establish that also $\mathcal{G}$ is a continuous image
of some $g(\mathcal{H}_m)$? In dimensions $0$ and $1$ it is. \\

{\bf D)} The Green function matrix $g_m(x,y)$ is still a bit of an enigma in the case when it is defined 
by a general function $h_m(x,y)$ defining the energy. We see for example that $g_2$ is identically zero
if $h_2$ is anti-symmetric. We also see that if $h(x,y)=1$ everywhere, then $g$ is invertible and positive
definite. 

\bibliographystyle{plain}

\end{document}